\documentclass[a4paper,11pt]{amsart}
\usepackage{etex}

\usepackage{verbatim}

\usepackage{multicol}
\usepackage[all]{xy}
\usepackage{appendix}
\usepackage{graphicx}
\usepackage{tikz-cd}
\usepackage{wrapfig}
\usepackage[pagebackref]{hyperref}
\hypersetup{
	colorlinks=true,
	linkcolor=cyan,
	citecolor=cyan
			}
\usepackage[makeindex]{imakeidx}

\usepackage{amsfonts, amsmath, amssymb, amsthm, amsopn, latexsym,
 graphicx, mathrsfs,  verbatim, enumerate, url, stmaryrd, enumitem, tikz,mathtools} 

\usepackage[normalem]{ulem}
\usepackage[all]{hypcap} 
\usepackage{pst-all} 
 \usepackage{lineno} 
 
 \usepackage{booktabs} 
\usepackage{xcolor}

\definecolor{green}{rgb}{0.0, 0.5, 0.0}
\definecolor{yellow}{rgb}{0.8, 0.33, 0.0}
\definecolor{purple}{RGB}{172, 51, 255}

\usepackage{caption}
\numberwithin{equation}{section}

\theoremstyle{plain} 
\newtheorem{theorem}[equation]{Theorem} 
\newtheorem{theorem*}{Theorem}

\newtheorem{lema}[equation]{Lemma}
\newtheorem{prop}[equation]{Proposition}
\newtheorem{cor}[equation]{Corollary}
\newtheorem{algorithm}[equation]{Algorithm}

\newtheorem{thmno}{Theorem}
 
\theoremstyle{definition} 
\newtheorem{defi}[equation]{Definition}

\newtheorem{remark}[equation]{Remark}
\newtheorem{example}[equation]{Example}
\newtheorem{notation}[equation]{Notation}
\newtheorem*{ack}{Acknowledgement}

\makeindex

\def\C {\mathbb C}
\def\N {\mathbb N}
\def\R {\mathbb R}

\def\Q {\mathbb Q}

\def\Z {\mathbb Z}

\def\a{\alpha}

\def\F{\mathcal{F}}

\def\fan{\mathcal{F}}

\def\l{\lambda}

\def\Newton{\mathcal N}

\def\M{\mathcal{M}}

\def\ord{\mathrm{ord}}

\def\r{\rho}
\def\reg{\mathrm{reg}}
\def\s{\sigma}
\newcommand{\slp}{\mathbf{S}} 

\def\prim{\mathrm{prim}}

\makeindex
\def\elem(#1,#2){  \left\{ \frac{#1}  {\overline {\ #2\ }} \right\} }  

\title{Multiplier ideals of plane curve singularities via Newton polygons}
\author[P.D. Gonz\'alez P\'erez]{Pedro D. Gonz\'alez P\'erez} 
\address{Instituto de Matematica Interdisciplinar y Departamento de \'Algebra, Geometr\'{\i}a y Topolog\'{\i}a, Facultad de Ciencias Matem\'aticas, Universidad Complutense de Madrid, Plaza de las Ciencias 3, Madrid 28040, Espa\~na.}
\email{pgonzalez@mat.ucm.es}
\author[M. Gonz\'alez Villa]{Manuel Gonz\'alez Villa}
\address{
Centro de investigaci\'on en 
Matem\'aticas\\
Apartado Postal
402 \\
36000 Guanajuato, Gto., M\'exico.
}
\email{manuel.gonzalez@cimat.mx}
\author[C.R. Guzm\'an Dur\'an]{Carlos R. Guzm\'an Dur\'an}
\email{guzman@cimat.mx}

\author[M. Robredo Buces] {Miguel Robredo Buces}
\address{Instituto de Ciencias Matem\'aticas (CSIC-UAM-UC3M-UCM), Calle Nicol\'as Cabrera,  13-15
Campus de Cantoblanco, 
28049 Madrid Spain}
\email{miguel.robredo@icmat.es}
\thanks{
The first and second author are supported by the grant PID2020-114750GB-C32.
The fourth author was supported by SEV-2015-0554}

\date{\today}


 
\begin{document}

\maketitle

\begin{quote}
\textbf{Abstract}
We give  a description of the multiplier ideals and jumping numbers associated with a plane curve singularity  in a smooth surface in terms of Newton polygons. 
Our approach is inspired by a theorem of Howald about multiplier ideals of Newton non-degenerate hypersurfaces and our results provide a generalization of it to the case of plane curve singularities.
We use toroidal embedded resolutions, which 
can be applied to the case of quasi-ordinary hypersurface singularities.
\end{quote}

\section*{Introduction}
Let $S$ be a smooth complex algebraic variety, $C$ be an integral effective divisor on $S$, and $\xi >0$ be a rational number. The \textit{multiplier ideal} associated with $C$ and $\xi$  is defined by
\[
	\mathcal{J}(\xi C) := \pi_{*}\mathcal{O}_{Y}(K_{\pi}-\lfloor  \xi \,  \pi^* (C) \rfloor ),
\]
where $\pi: Y \to S$ is a log-resolution of $C$ and $K_{\pi}$ is the relative canonical divisor. 
For any point $o \in S$ there  exists an increasing sequence $(\xi_{i})_i$ of positive rational numbers, called 
the \textit{jumping numbers} of $C$ at $o$, 
such that if  $\xi_i \leq \xi <   \xi_{i+1}$ then
$
\mathcal{J} (\xi C )_o = 
	\mathcal{J} (\xi_{i} C )_o
	\supsetneq \mathcal{J}( \xi_{i+1} C)_o
$. 
More generally, the multiplier ideals and their jumping numbers can be associated 
with an ideal sheaf on $S$, even  if $S$ has some mild singularities. 
The study of multiplier ideals have become a central aspect of birational geometry, thanks to 
the vanishing theorems of Kawamata, Viehweg and Nadel, which were inspired by 
Kodaira vanishing theorem. 
We refer to Lazarsfeld's book \cite{PAG} for the historical aspects and more information on the subject.
Many properties of the jumping numbers were studied by
Ein, Lazarsfeld, Smith and Varolin \cite{ELSVJumpingCoefficients}.
The multiplier ideals provide a subtle measure of the singularities of the pair $( S, C)$ and enjoy 
a wealth of relations with other notions 
 like the Hodge modules,  mixed Hodge modules, 
 $V$-filtrations of $D$-modules 
 and poles of Igusa zeta functions which are of interest in singularity theory 
 (see Budur's survey \cite{BudurSingulInvMilnorFiber}). 
 The multiplier ideals can be defined in complex analytic terms without using any resolution of singularities. 
Favre and Jonsson studied multiplier ideals associated 
with an ideal of the ring of germs of holomorphic functions at a point of 
a smooth surface by using \textit{tree potentials} on the \textit{valuative tree}
(see \cite{FJ05, TVT}  and Jonsson's survey  \cite[Section 7]{JonssonDynamicsBerkovich}).
Multiplier ideals are also connected with the study of singularities over fields of positive characteristic through 
test ideals (see for instance \cite[Section 4]{SchTuck}). 
Algorithms to  compute multiplier ideals and jumping numbers were given Shibuta
 \cite{Shibuta11}, and  Berkesch and Leykin \cite{Berkesch10}. 
 
 \medskip 
If $S$ is a complex surface there are many results about multiplier ideals and
 their associated jumping numbers.
J\"arvilehto described 
the jumping numbers of a simple complete ideal in a two-dimensional regular local ring, and 
deduced a formula for the jumping numbers of a \textit{branch} $(C,o)$ on a smooth surface $S$ (see \cite{JarJNSCI}). Smith and Thompson introduced the notion of \textit{contribution} of an 
exceptional prime divisor of 
the minimal embedded resolution $\pi$ of a plane curve singularity $C$
on a smooth surface $S$, and they relate this notion 
with the rupture components of the divisor $\pi^*(C)$ (see \cite{STIrrelExcDiv}).  
Tucker \cite{TuckerJNAlgSurf} and 
Naie \cite{JNN} gave also descriptions of the jumping numbers in this case. 
Tucker studied multiplier ideals on a surface $S$ with rational singularities by using 
the notion of 
\textit{critical contribution} of a \textit{reduced} exceptional divisor to a jumping number and 
applied it to give an algorithm to compute the jumping numbers (see \cite{TuckerJNRatSing, TuckerJNAlgSurf}). 
Moreover, he established an iterative relation between the jumping numbers of a branch and those of its
approximate roots (see   \cite[Chapter VI]{TuckerJNAlgSurf}).
Alberich-Carrami\~nana,  \`Alvarez Montaner, and Dachs-Cadefau gave a different algorithm 
to compute the jumping numbers and the associated multiplier ideals (see \cite{AAD}). 
Later on Alberich-Carrami\~nana, \`Alvarez Montaner and Blanco
 gave an algorithm to compute the integral  closure 
of an ideal of $\C \{ x, y \}$
in terms of monomials in a set of maximal contact elements  
of  its minimal log-resolution. Combining this algorithm with the one in
\cite{AAD},  they obtained an algorithm 
to compute the multiplier ideals and jumping numbers of an ideal  of $\C \{ x, y \}$ in
terms of a system of generators of it (see \cite{Blanco}). 
In the case of an analytically irreducible plane curve singularity, 
a similar result describing the associated multiplier ideals 
was obtained with different methods in \cite{GuzmanPhD} and also in Zhang's preprint \cite{Zhang}. 
Other recent works about jumping numbers of multiplier ideals in this context are   
 \cite{HJFormulaJN, HJHNORdTree}.

\medskip 

The aim of this paper is to develop an algorithmic and conceptually new description of the multiplier ideals associated with a plane curve singularity on a smooth surface. We give a combinatorial algorithm to provide the generators of the multiplier ideals and the jumping numbers associated with a plane curve singularity on a smooth surface in terms of a finite set of Newton polygons 
appearing in a  toroidal embedded resolution process of the plane curve singularity. Our approach is inspired by Howald's result  about multiplier ideals of hypersurface singularities which are Newton non-degenerate (see \cite{Howald:Non-degenerate}), but we 
do not require any Newton non-degeneracy hypothesis.
Our motivation was to develop a method which can be extended to the study of higher-dimensional singularities, as irreducible quasi-ordinary hypersurface singularities (see  \cite[Ch. 6]{PhD-Robredo}). 
In contrast to previous results in the literature, our techniques do not pass through the correspondence between antinef divisors and complete ideals, used in  \cite{Blanco} for instance. This paper is a development of the PhD Theses of the
third  and fourth named authors \cite{GuzmanPhD,PhD-Robredo}.

\medskip

In order to state our results, we  outline briefly the construction of 
a toroidal embedded resolution of a plane curve singularity  $(C,o)$ embedded in a complex smooth surface $ (S, o)$
These kinds of resolutions are described in Subsection \ref{subsec:torres},
following the presentation given by 
Garc\'\i a Barroso, Popescu-Pampu and the first named author in  \cite{GBGPPP20}. The process of toroidal embedded resolution that we recall here 
is a slight generalization of
the toroidal resolution processes of plane curves  by  
Oka \cite{OkaToroidalRes}, L\^e and  Oka \cite{ResolComplex}, A'Campo  and Oka \cite{ACO96}, the first named author \cite{GPRQo} 
and, Cassou-Nogu\'es 
and Libgober \cite{CNL14}.

\medskip

\medskip
Let us denote by $\mathcal{O} \cong \C \{ x , y \}$ the local ring of germs of holomorphic functions of 
$S$ at $o$. 
We start by fixing a cross $(R, L)$ at $o$, which is an ordered pair of smooth transversal branches 
defined by the vanishing locus of the entries of a local coordinate system $(x, y)$ of $S$ at $o$. We assume that $R$ is not a component of $C$.
We denote 
by $\Newton_{R, L} (C)$ 
the Newton 
polygon of any defining function $f_C \in \mathcal{O}$ of $C$. 
Then, we consider the \textit{regularized Newton modification} of $C$ with respect to the cross $(R, L)$, 
which is a toric modification defined in terms of the Newton polygon $\Newton_{R, L} (C)$. 
Its exceptional curve $E$ intersects the strict transform of $C$ only at smooth points of $E$. If 
$o'$ is any of these points, we denote by $R'$ the germ of $E$ at $o'$, and 
we choose a cross $(R', L')$ at $o'$ by taking $L'$ as a  \textit{curvetta} of $E$ at $o'$.
Then, we iterate this procedure until we get a toroidal embedded resolution $\pi: \Sigma \to S$  of $C$. 
In the process, we have a finite number of crosses $(R_i, L_i)$ 
at some infinitely near points $o_i$ of $o$, for $i \in I_\pi$.
The completion $\hat{C}_\pi$ of $C$ relative to $\pi$
is the reduced plane curve whose branches are $R$ and the projections by $\pi$ of the smooth branches $L_i$, for $i \in I_\pi$.
 By definition every branch of  $C$ is a branch of $\hat{C}_\pi$.  
 
It is possible to choose $\pi$ in such a way
that $\pi$ is the minimal embedded resolution of $C$ and 
the components of $\hat{C}_\pi$ are 
components of $C$ or maximal contact curves of $\pi$ (see Definition  \ref{mcc} and Remark \ref{rem:red}).

Let us fix any toroidal embedded resolution $\pi$ of $C$.
If $h \in \mathcal{O}$
defines a plane curve singularity $C_h \subset S$,  we denote by 
$\Newton_{R_i, L_i} (C_h)$ 
the Newton polygon of the total transform of $C_h$ with respect to 
the cross $(R_i, L_i)$ and by 
$\lambda_{R_i}$  the log-discrepancy of 
$R_i$, for $i \in I_\pi$. In Theorem \ref{Thm:Howalditerated} 
we prove that: 
\begin{thmno} 
The multiplier ideal 
$\mathcal{J}(\xi C)_o$ consists of the functions $h \in \mathcal{O}$ such that   the inclusion
\[	\Newton_{R_i, L_i} (C_h) + (\lambda_{R_i}, 1)  \subset \mathrm{Int} (\xi \Newton_{R_i, L_i} (C)) \]
is satisfied for all $i \in I_\pi$.
\end{thmno}

The proof is based on a basic property of Newton polyhedra, which is combined 
with the description of the log-discrepancies of the exceptional divisors in the toroidal embedded resolution 
(see Proposition \ref{mcp}). 
If $C$ is Newton non-degenerate with respect to  the coordinate system $(x, y)$
then we recover a particular case of Howald's Theorem \cite{Howald:Non-degenerate}.

\medskip

By convenience, if  $D$ is an exceptional prime divisor (resp. $D$ is a branch on $S$) 
we denote by $\nu_D$ the associated divisorial valuation (resp. vanishing order valuation). 
We prove that  a function $h \in \mathcal{O}$ belongs to the multiplier ideal $\mathcal{J}(\xi C)_o$
if and only if 
$\nu_{D} (C_h)  + \l_D > \xi \nu_{D}(C)$,
for $D$ running through the 
rupture components of the divisor $\pi^*(C)$ and 
the components of $C$ (see Corollary  \ref{prop-form}). 
Denote by $\{x_0, \dots, x_s \}$ the set obtained by taking a defining function for every branch of $\hat{C}_\pi$.
Every monomial $\M$ in $x_0, \dots, x_s$ determines the rational number
\[
\xi_{\mathcal{M}}:= 
	\min_{D} = 
		\{ ({\nu_{D}(\mathcal{M})+\lambda_{D}}) ({\nu_{D}(C)} )^{-1}  \}, 
\]
where $D$ runs through the rupture components of the divisor $\pi^*(C)$ 
and the branches of $C$.

\medskip 

The following theorem 
is based upon previous works of Spivakovsky \cite{VS}, 
Delgado, Galindo, and Nu\~nez \cite{DelGalNuGenSeq}
and Robredo \cite{PhD-Robredo} about generating sequences of valuations (see Theorem \ref{Thm:Monomiality}). 
\begin{thmno} 
Every jumping number of $C$ is of the form 
$\xi_\M$ and 
the multiplier ideal 
$\mathcal{J}(\xi C)_o$ is generated by the finite set of  monomials 
$\M$ in $x_0, \dots, x_s$ such that 
$\xi < \xi_{\M} \leq \xi +1$.
\end{thmno}

As a consequence of our main results, the computation 
of the multiplier ideals and the associated jumping numbers of $C$ 
boils down to an optimization problem in terms of log discrepancies 
of the rupture components of $\pi^*(C)$ and the values of the corresponding exceptional divisors on the functions $x_0, \dots, x_s$.

\medskip 
We have formulated our results in the complex analytic category, but they 
 also hold for algebroid curves on a smooth surface over algebraically closed fields of 
arbitrary characteristic (see Remark \ref{rem:char}). 

\medskip

The structure of the paper is as follows. 
In Section \ref{Sect:Definition mult ideals and jumping numbers} we recall the basic notions 
about multiplier ideals. In Section \ref{Chapter:General background}
 we introduce some notations and well-known results about 
plane curve singularities and valuations. 
In Section \ref{sec:tor-res} we describe an algorithm of 
toroidal embedded resolution following \cite{GBGPPP20}.
The main results of the paper are proven in Section \ref{sec:main}, where we  illustrate our results with a detailed example. 

\section{Multiplier ideals and jumping numbers}
\label{Sect:Definition mult ideals and jumping numbers}

In this section we briefly review basic definitions and properties of the theory of multiplier ideals. For further details we refer to \cite[Chapter 9]{PAG}.

Let  $X$ be a 
 smooth complex algebraic variety and let  $\mathfrak{a}\subseteq \mathcal{O}_{X}$ be an ideal sheaf. 
A \textbf{log-resolution} of $\mathfrak{a}$ is a modification (proper and birational map)
$\pi:Y\rightarrow X$, with $Y$ smooth, 
exceptional locus $E$, and such that
$\pi^{*}\mathfrak{a}=\mathcal{O}_{Y}(-F)$,
where $F$ is an effective divisor such that 
$F+E$ has simple normal crossings. 

\medskip

For  $a \in \Q$, we denote by $\lfloor a \rfloor$ the greatest integer lower than or equal to $a$. 
For a $\Q$-divisor  $D = \sum_j {a_j } D_j$,  supported on the prime divisors $D_j$, we denote by $\lfloor D \rfloor \coloneqq \sum_j \lfloor {a_j }  \rfloor D_j$.

\medskip 
Let $\pi :Y\rightarrow X$ be a log-resolution of an ideal sheaf $\mathfrak{a}$ of $\mathcal{O}_X$. Denote by $\boxed{K_{\pi}}$ the relative canonical divisor, which is equal to the divisor 
associated with the jacobian determinant of $\pi$.
The \textbf{multiplier ideal} sheaf $\boxed{\mathcal{J}(\mathfrak{a}^{\xi})}$
associated to $\xi \in\mathbb{Q}_{> 0}$  and $\mathfrak{a}$ is defined as $
\mathcal{J}(\mathfrak{a}^{\xi}) \coloneqq 
	\pi_{*}\mathcal{O}_{Y}(K_{\pi}-\lfloor  \xi F\rfloor )$.

 \medskip 

The definition of the multiplier ideal $\mathcal{J}(\mathfrak{a}^{\xi})$ relies on the choice of a 
log-resolution of $\mathfrak{a}$, but  it is independent of it (see \cite[Theorem 9.2.8]{PAG}).
The multiplier ideal $\mathcal{J}(\mathfrak{a}^{\xi})$ can be characterized  in terms of valuations. If $E_i$ is a prime divisor on $Y$ we denote by 
$\nu_{E_i}$ the vanishing order valuation along $E_i$. 
A prime divisor $E_i$  contained in the support of $F + K_\pi$ is either the strict transform of a divisor on $X$ or  must be contained in the exceptional divisor of $\pi$.
Let us write 
$F		=  \sum r_{i}E_{i}$,
and
$ K_{\pi} = \sum (\lambda_{E_{i}}-1) E_{i}$,  
where the $E_{i}$ are the prime divisors in the support of $E+F$ on $Y$. 
Then, 
\begin{align} \label{Jac}
\mathcal{J} (\mathfrak{a}^\xi) = 
	\{ h \in \mathcal{O}_{X} \mid  \nu_{E_i}(h) \geq \lfloor \xi r_{i} \rfloor - (\lambda_{E_{i}}-1) 
		\mbox{ for all } i\}, 
\end{align}
or, equivalently,
\begin{align} \label{eq:form MI}
\mathcal{J} (\mathfrak{a}^ \xi )  = 
	\{ h \in \mathcal{O}_{X} \vert \mbox{ }
	 \nu_{E_i}(h) + \lambda_{E_{i}} > \xi \mbox{ }  r_i 
		\mbox{ for all } i \} .  
\end{align}
The equivalence follows  since for any $a \in \Z$ and $b \in \Q$ it holds that $a \geq \lfloor b \rfloor$ if and only if $a > b-1$.
The number $\lambda_{{E}_i}$ is called the \textbf{log-discrepancy} of the exceptional prime divisor $E_i$. 
 If $B \subset Y$ is the strict transform of a prime divisor in $X$ we set $\lambda_{B} :=1$.

\medskip 

 The next lemma introduces some numerical invariants of set of  multiplier ideals of an ideal sheaf $\mathfrak{a}$ of $\mathcal{O}_X$ at a point $x \in X$. 
\begin{lema}[See  \cite{PAG}, Lemma 9.3.21]  
\label{Def:Jumping numbers}  
Let $X$ be an smooth algebraic variety and let $\mathfrak{a} \subseteq \mathcal{O}_{X}$ be an ideal sheaf and 
$x \in X$.
There exists a strictly increasing discrete sequence $(\xi_{i})$ of positive rational numbers 
such that if  $\xi \in \Q \cap [\xi_{i}, \xi_{i+1})$, then $
\mathcal{J}(\mathfrak{a}^{\xi_{i}})_x = 
	\mathcal{J}(\mathfrak{a}^{\xi})_x \supsetneq \mathcal{J}(\mathfrak{a}^{\xi_{i+1}})_x $. 
\end{lema}

\medskip 
The numbers $\xi_{i}$ are called the \textbf{jumping numbers} associated with $\mathfrak{a}$ at $x$. The smallest jumping number $\xi_{1}$ is called the \textbf{log-canonical threshold} of $\mathfrak{a}$ at $x$.

\medskip 
If $D$ is an effective integral divisor on $X$, determining the line bundle $\mathcal{O}_X (-D) =
	\{ h \in \mathcal{O}_X \mid \mathrm{div} (h) - D \geq 0 \}$, 
	we denote by $\boxed{\mathcal{J}(\xi D)}$ 
	the  multiplier ideal 	$\mathcal{J} (\mathcal{O}_X (-D)^\xi )$ 
associated with the ideal sheaf $\mathcal{O}_X (-D)$ and the number $\xi$.
The following periodicity property of  the multiplier ideals  $\mathcal{J}(\xi D)$  imply that  their jumping numbers are determined by the finitely many of them lying in the unit interval $(0,1]$
(see   \cite[Example 1.7 and Remark 1.15]{ELSVJumpingCoefficients}  or \cite[Example 9.2.12 and Proposition 9.2.31]{PAG}).
\begin{lema} \label{Lem:periodicidad} 
Let $D$ be an effective integral divisor on a smooth variety $X$.  
		Then, $\mathcal{J} \left(  D \right) =  \mathcal{O}_X (-D) $,
  the number  $1$ is a jumping number of the multiplier ideals of $D$ and
		$\mathcal{J} \left( (\xi + 1 ) D \right) = \mathcal{J} \left( \xi  D \right) \otimes \mathcal{O}_X (-D)$.		
\end{lema}

\begin{remark} 
If $X$ is a smooth variety and $f$ is a germ of complex analytic function at $x \in X$,
then the above definitions of multiplier ideals and jumping numbers of $(f)$ 
generalize to this local setting (see \cite[Remark 1.26]{ELSVJumpingCoefficients}). 
\end{remark}

\section{Basic notions about plane curve singularities} \label{Chapter:General background}

Let $\boxed{S}$ be a smooth complex algebraic or analytic surface and $o \in S$ a fixed closed point. 
Denote by   $\boxed{\mathcal{O}}$, 
the local ring of germs of holomorphic functions on $S$ at $o$.

\medskip 

A germ of a complex analytic curve $\boxed{(C,o)}$ is defined by an equation,  $f_C = 0$,
where $\boxed{f_C} \in \mathcal{O}$ is a representative of $(C,o)$
(defined up to multiplication by a unit of the ring $ \mathcal{O}$). 
Similarly, if $f \in  \mathcal{O}$ we denote  $\boxed{C_f}$ the curve germ
defined by $f=0$.
If the base point $o$ of the germ is clear from the context we denote $(C, o)$ simply by $C$. 
The curve $C$ is a \textbf{branch} if $f_C$ is irreducible. 
 In general, if we decompose  $f_C = f_{1}^{a_1} \cdots f_{r}^{a_r}$ as a product of irreducible elements in
$\mathcal{O}$ with $f_i$ non-associated to $f_j$ for $i \ne j$,  and we put $C_i = C_{f_i}$ 
then we represent
$C = \sum_{i=1}^r a_i C_{i}$
as the effective divisor defined by $f_C$.

\medskip 

A \textbf{cross} at the point $o$ of $S$ is an ordered pair $(R, L)$ of transversal smooth branches. 
A local coordinate system $(x,y)$ on $(S,o)$ is an ordered pair of elements generating the maximal ideal of $\mathcal{O}$.
It defines a cross $(R, L)$ with $R=C_x$, $L=C_y$.

\medskip

A \textbf{model} of $(S,o)$ is 
a proper and birational map $\pi: (S_\pi, E_\pi) \to (S,o)$
such that the restriction of $\pi$ to $S_\pi \setminus   E_\pi
\to S \setminus \{ o \}$ is an isomorphism, where $E_\pi = \pi^{-1} (0)$. 
If $\pi$ is not the identity map, it is a composition of blow ups of points infinitely near to $o$, 
hence the irreducible components the reduced divisor $\boxed{E_\pi} := \pi^{-1} (o) $ are
projective lines.   We denote by $\boxed{E(\pi)}$ the  set of prime divisors of the exceptional divisor $E_\pi$, that is, we can write
$E_\pi = \bigcup_{i \in E(\pi)} E_i$.

\medskip 

\begin{defi}
Let $\pi$ be a model of $(S,o)$ and let  $C$ be a plane curve germ on $S$ at  a point $o$. The  \textbf{total transform} $\pi^* ( C)$ is the divisor of  $f_C \circ \pi$. The \textbf{strict transform} $\boxed{C^\pi}$ of $C$ by $\pi$
is sum of components of $\pi^* ( C)$ which are supported on 
the closure of $\pi^{-1} (C \setminus \{ o \})$. 
The model $\pi$ is an
\textbf{embedded resolution} of $C$
(also called 
\textbf{log-resolution}) if $\pi^* ( C)$ has simple normal crossings.

\end{defi}

 Any embedded resolution of $C$ is a composition of blow ups of a finite set of infinitely near points of $o$. There exists a unique \textbf{minimal embedded resolution} of $C$, which is the one requiring the smallest number of blow ups.

\medskip

We introduce several notions of dual graph associated to 
certain divisors on a smooth surface. 
\index{dual graph} \index{valency of a vertex} \index{rupture divisor} 
 Let $D = \sum_{j \in J} D_j$ be a reduced divisor with  simple normal crossings on a
smooth surface $\Sigma$.
The \textbf{dual graph} $\boxed{G (D)}$  \index{dual graph} of $D$ is
the combinatorial graph  with vertex set $J$ and whose edges  
are in bijection with 
the singular points of $D$. If $p$ is a singular point of $D$ then 
there are unique elements $j, k \in J$ such that $p \in E_j \cap E_k$. 
Then $\mathcal{E}_p = \{j, k\}$ is the corresponding edge of $G(D)$. 
The \textbf{valency} \index{valency}
of a vertex $j$ of $G (D)$ is the number of edges 
incident to $j$. 
If $D$ is an effective divisor
we denote also by $\boxed{G (D)}$
the dual graph associated to the reduction of $D$.
A prime divisor $D_j$ whose corresponding vertex in the dual graph  $G (D)$ 
has valency $\geq 3$ is called a \textbf{rupture component} of $D$.
 We denote by $\boxed{\mathcal{E} (D)} $ the set consisting of 
prime divisors of $D$ defining end vertices of the dual graph $G(D)$, that is, vertices of valency one.

\begin{defi} \label{def:DG}  
Let $\pi$ be a model of $(S,o)$. 
The dual graph \index{dual graph} $\boxed{G(\pi)}$
of $\pi$ is  $G (E_\pi) $. 
If $\pi$ is an embedded resolution of a plane curve $C$, 
we denote by $\boxed{G(\pi, C)}$ the dual graph of  $\pi^*(C)$,
 by $\boxed{\mathcal{R}_\pi (C)}$ the set of \textbf{rupture components} of the divisor $\pi^*(C)$.  
\end{defi}

The dual graphs $G( \pi, C)$ and $G( \pi)$ are finite trees, that is, 
connected graphs with a finite number of vertices and with no cycles.

\medskip 

Let $\pi$ be a model of $(S, o)$.
A branch $K_i$ in $S$ is a \textbf{curvetta} at a component $E_i$ of the exceptional divisor $E_\pi$ 
if $K_i^\pi + E_\pi $ is a simple normal crossing divisor such that 
$K_i ^\pi \cap E_i \ne \emptyset$.
In particular, in this case $\pi$ is an embedded resolution of $K_i$.
We define now some classes of finite subsets whose elements are curvettas.

\begin{defi} \label{mcc} \index{set of maximal contact curves}
Let $\pi$ be a model of $(S,o)$, different from the identity map of $S$
or the blow up of $o$ and let $C$ be a plane curve germ on $S$ at $o$.
\begin{itemize}
\item 
A set of \textbf{maximal contact curves} of $\pi$ contains exactly one 
 curvetta 
$K_i$ at $E_i$, for every $E_i \in \mathcal{E} (E_\pi)$. 

\item If $\pi$ is an embedded resolution of $C$, a set of \textbf{maximal contact curves} of the pair $(\pi, C)$
consist of  the components of $C$ 
together with one curvetta 
$K_i$ at $E_i$, for $E_i$ running through the components of $E_\pi$ 
in the set $\mathcal{E} (\pi^*(C))$.
\end{itemize}
\end{defi}

Let us recall some facts about the valuations which we use in this paper.

\begin{defi} \label{semivaluation} \index{semivaluation} \index{valuation}
A {\bf valuation} of the local ring $\mathcal{O}$  is a function 
$\boxed{\nu}: \mathcal{O} \rightarrow [0, \infty]$ such that
\begin{enumerate}  [label=({\alph*})]
\item \label{add} $\nu(fg)=\nu(f) + \nu(g)$, 
  \item \label{ineq} $\nu(f+g) \geq \min(\nu(f), \nu(g))$ for all $f,g \in \mathcal{O}$;
\item \label{const} 
$\nu (f) = \infty \Leftrightarrow f =0$, 
\end{enumerate} 
where $ [0, \infty] := \R_{\geq 0} \cup \{ \infty \}$ is considered
as a semigroup with addition in the usual sense.
\end{defi}
 If $D$ is a germ of curve on $S$ 
and $\nu$ is a valuation of $\mathcal{O}$ then we denote $\boxed{\nu (D)}:= \nu (f_D)$.

\medskip 

Let us introduce two useful types of valuations associated to a branch $C$ or to an exceptional prime divisor.

\medskip 
The \textbf{vanishing order valuation  along a branch $C$}, denoted by $\boxed{\nu_{C}}$, is given by $\nu_{C} (h) = a$, if  $0 \ne h \in \mathcal{O} $, and $h=f_C ^{a}g$ with  $\mathrm{gcd} (f_C, g) = 1$.  We set $\nu_{C} (0) = \infty$. 
Notice that the vanishing order valuation along a branch  is well defined because  $\mathcal{O}$ is a unique factorization domain.

\medskip 
Let  $\pi$ be a model of $S$ at $o$ and $E_i$ be an irreducible component of $E_\pi$.  The function 
$\boxed{\nu_{E_i}}: \mathcal{O} \to [0, \infty]$, which maps $h \in \mathcal{O}$ to 
 the order of vanishing along $E_i$ of
$h \circ \pi$, is a valuation 
of $\mathcal{O}$, called the \textbf{divisorial valuation} of $E_i$.

\medskip

Denote by $\boxed{\succeq}$ the poset relation in $(\R \cup \{ \infty \})^s$ given by $ (a_1, 
\dots, a_s)   \succeq  (b_1, \dots, b_s) $ if 
$a_{i} \geq b_{i}$ for every index $i \in \{ 1 , \ldots , s \}$.
Let $\nu_{1},\ldots ,\nu_{s}$ be  valuations of $\mathcal{O}$. 
The \textbf{valuation ideal} $\boxed{I_{\underline{c}}^{\underline{\nu}}}$ associated with $\underline{c} =(c_1, \dots, c_{s}) \in \R^{s}_{\geq 0}$ and $\underline{\nu}:= (\nu_{1},\ldots ,\nu_{s}) $ is 
$ {I_{\underline{c}}^{\underline{\nu}}}:= \{ f \in \mathcal{O}  \mbox{ } \mid \mbox{ }  	\underline{\nu} (f) \succeq\underline{c} \}$.
\medskip 

The notion of generating sequence was considered in \cite[Definition 1.1]{VS} for one valuation, and studied in \cite{CG03, DelGalNuGenSeq} for tuples of divisorial valuations. 
\begin{defi} \label{Def:Generating sequence valuation tuple}
Let  $\underline{\nu} = (\nu_{1},\ldots ,\nu_{s})$ be a tuple of valuations of $\mathcal{O}$. A set of elements $\{ z_{j} \}_{j \in J}$  in the maximal ideal of $\mathcal{O}$
is a \textbf{generating sequence} \index{generating sequence}
of $\underline{\nu}$ if for every $\underline{c} \in \R^s_{\geq 0}$
the valuation ideal $I_{\underline{c}}^{\underline{\nu}}$ is generated by the set
$
\{ {\underset{{{\tiny \hbox{finite}}}}{\prod}} z_{j}^{b_{j}} \mbox{ } \mid   
	b_{j} \in \Z_{\geq  0}, \mbox{ } \sum b_{j} \,  \underline{\nu} (z_{j})  \succeq  \underline{c} \}.
$	
A generating sequence is \textbf{minimal} if every proper subset of it fails to be a generating sequence. 

For instance, a minimal generating sequence of the divisorial valuation associated with the exceptional divisor of the blow up of $o$ in $S$ is 
$(x, y)$, where $(x, y)$ is a local system of coordinates of $S$ at $o$.

\end{defi}

Let  $\underline{\nu} = (\nu_{1},\ldots ,\nu_{s})$ be a tuple of divisorial valuations. Campillo and Galindo proved that $\underline{\nu}$ has 
a finite generating sequence (see \cite[Th. 3]{CG03}). 
If $\{ z_{j} \}_{j \in J}$ is a finite generating sequence of $\underline{\nu}$, 
then it is a generating sequence of $\nu_1$ (see \cite{DelGalNuGenSeq}), in
particular, it is a set of generators of the maximal ideal of $\mathcal{O}$
(see \cite{VS}).
A model $\pi:  (\Sigma, E) \to (S, o)$ is a \textbf{minimal embedded resolution} 
of $\underline{\nu}$ if $\pi$ is a composition of the minimal number of blowing ups of 
points such that  there exists a component $E_i$ of the exceptional divisor $E$ such that  $\nu_i = \nu_{E_i}$, for $i= 1, \dots,s$. 
 We have the following characterization of the minimal generating sequences associated
with a tuple of divisorial valuations:

\begin{theorem}\cite[Th. 5]{DelGalNuGenSeq} \label{Thm:DGN} 
Let $\underline{\nu} = (\nu_{1},\ldots ,\nu_{s})$ be a tuple of divisorial valuations 
of $\mathcal{O}$. We assume that if $s =1$ then $\nu_1$ is different from divisorial valuation associated with the blow up of $o$ in $S$.  Let us denote by $\pi:\Sigma \to S $ the
minimal embedded resolution 
of $\underline{\nu}$. 
Take a defining function for every curvetta in a set of maximal contact curves of $\pi$. Then, these functions define a minimal generating sequence of $\underline{\nu}$.
\end{theorem}

\begin{cor}   \label{Cor:DGN}
Let $C = \sum_{j=1}^r a_j C_j $ 
be a singular plane curve germ on $S$ at $o$. 
Let $\pi: \Sigma \to S$ be the minimal embedded resolution of $C$. 
Take a defining function for every element in a set of maximal contact curves of $\pi$  (resp. of the pair $(\pi, C)$). Then, these functions define a minimal generating sequence of the set of divisorial 
valuations associated with 
 the rupture components of the divisor  $\pi^*(C)$ (resp. of the set of divisorial 
valuations associated with 
the rupture components of the divisor $\pi^*(C)$
and the vanishing order valuations of
the branches of $C$). 
\end{cor}
\begin{proof} 
 We have that $\pi$ is also the minimal resolution of the divisorial valuations of 
the rupture components of the divisor $\pi^*(C)$. Let us denote by 
$\tilde{C}$ the sum of maximal contact curves of $(\pi, C)$, considered as a reduced effective divisor.
We have that 
$\pi$ is also the minimal embedded resolution of $\tilde{C}$.
The dual graph $G(\pi)$ is obtained from 
$G(\pi, \tilde{C})$ by 
deleting the arrows corresponding to the components 
of $\tilde{C}$. 
The first assertion follows by  Theorem \ref{Thm:DGN}, since 
a set of maximal contact elements
of $(\pi, C)$ contains 
a set of maximal contact elements of $\pi$. 
We refer to 
\cite[Cor. 4.160]{PhD-Robredo}
for a proof of the second statement. 
\end{proof}

\section{Toroidal embedded resolutions of plane curves} \label{sec:tor-res}

In this section we recall the construction of a toroidal embedded resolution 
of a plane curve germ following  an algorithm given in \cite{GBGPPP20}. 
We start by fixing some vocabulary about toric geometry 
and Newton polygons in the two dimensional situation.
We describe the properties of the Newton modification associated to a plane curve singularity and a cross. 
Then we recall 
the construction of {\green a} toroidal embedded resolution 
of a plane curve germ and its associated 
combinatorics which is encoded in the associated fan tree.

\subsection{Fans, cones and toric varieties} \label{sec:fan}

We introduce some standard notations of toric geometry following  \cite{Cox,FTV,OdaCB,EwCombConv}.

\medskip

A lattice $\boxed{N}$   is a free group of finite rank $d$. We denote by $\boxed{N_\R} : = N \otimes_\Z \R$ the real vector space spanned by $N$ and 
by $M_\R$ and $M$ its duals respectively. We denote by 
\[
\boxed{\langle, \rangle} \colon N_\R \times M_\R \to \R, \quad (u, v) \mapsto \langle u, v \rangle = v(  u) 
\]
the duality pairing between these two vector spaces. 
A cone $\theta \subset N_\R$  is
\textbf{rational} with respect to $N$ if it is 
of the form $\theta = \R_{\geq 0} v_1 + \cdots + \R_{\geq 0} v_s$ for $v_1, \dots, v_s \in N$ and 
$s \in \N$.
It is \textbf{strictly convex} if $\{0\}$ is the biggest 
subspace contained in it. 
A face of the cone $\theta$ is the intersection of it 
with a  supporting hyperplane, that is a subspace of codimension one such that one of its
half-spaces contains 
$\theta$. 
A cone is \textbf{regular} if it is spanned  by a subset of a basis of the lattice $N$, 
where by convention the cone generated by the empty set is $\{0\}$. 
A \textbf{fan $\boxed{\mathcal{F}}$  of the lattice} $N$ is a finite set of strictly   convex rational cones in $N_\R$ 
such that it is closed under the operation of taking faces of its cones and the intersection of any two cones in the fan is 
a face of each of them. The \textbf{support}  $\boxed{| \mathcal{F} |}$ of the fan  $\mathcal{F}$ is the union of its faces. 
We denote by $\boxed{\F_{\prim}}$ the set of primitive integral vectors of 
the lattice $N$ which spans the rays of the fan $\F$. 
The fan $\F$ is \textbf{regular} if all its cones are regular. 

\medskip

If $\theta \in \F$ we denote by  $\check{\theta} = \{v \in  M_\R \mid \langle u, v \rangle \geq 0, \forall u \in \s \}$ the dual cone of $\theta$.
The coordinate ring  
of the affine toric variety $\boxed{X_\theta}$ is isomorphic to the semigroup algebra  
$\C [\check{\theta} \cap M] = \{ \sum c_\a \chi^\a \mid \a \in \check{\theta} \cap M , \, c_\a \in \C \}$. 
We denote by $\boxed{X_{\F}}$ the normal toric variety associated to the fan $\F$, 
which is obtained by 
glueing up the varieties $X_\theta$ for $\theta \in \F$. 
There is an action of the torus $\boxed{T_N}:= X_{\{ 0 \}}$ on $X_\F$, and 
the orbits of this action are in bijection with the cones in the fan $\F$. 
If $\theta \in \F$ we denote by $\boxed{O_\theta}$
the corresponding orbit of torus action on $X_\F$. 
The toric variety $X_{\F}$ is non-singular if and only if the fan $\F$ is regular. 
If $\theta \subset N_\R$ is a strictly convex cone which is rational for the lattice $N$, the set of faces of $\theta$ 
is a fan and the associated toric variety  coincides with $X_{\theta}$.

\medskip 
 If $\theta', \theta \subset N_\R$  are strictly convex rational cones  such that $\theta' \subset \theta$,
then we have a toric morphism
$
\psi^{\theta'}_\theta: X_{\theta'} \to X_{\theta}
$, which is 
determined by the inclusion of semigroups 
$\check{\theta} \cap M \to  \check{\theta'} \cap M$.

\medskip 

Let $\mathcal{F}$ and $\mathcal{F}'$ be two fans of the lattice $N$. The fan $\mathcal{F'}$ is 
\textbf{subdivision} of $\F$ if $|\F| = |\F'|$ and 
for every $\theta' \in \F'$ there exists $\theta \in \F$ such that $\theta' \subset \theta$. 
If $\F'$ is a subdivision of $\F$ we have a toric morphism
$\boxed{\psi^{\F'}_{\F}} : X_{\F'} \to X_{\F}$, obtained by glueing 
the maps $\psi^{\theta'}_\theta$ for every $\theta \in \F$, $\theta' \in \F'$
such that $\theta' \subset \theta$. The morphism $\psi^{\F'}_{\F}$  is 
a modification that is, it is proper and birational.
In addition, $\mathcal{F}'$ is a \textbf{regular subdivision} of $\F$ if 
$\F'$ is a regular fan containing every regular cone of $\F$.

\medskip 

We apply now these notions when the lattice $N$ has rank two and 
we endow it with a fixed basis $e_1, e_2$. 
We denote by  $\boxed{\s_0}$ the regular cone spanned by the basis $e_1, e_2$ of $N$. 
Note that $X_{\s_0} = \C^2$ and  $\C[X_{\s_0}] = \C[ x, y ]$ where $x = \chi^{\check{e}_1}$, $ y = \chi^{\check{e}_2}$ are the characters associated with the dual
basis  $\check{e}_1, \check{e}_2$ of $e_1, e_2$. 
If $\F$ is a fan of $N$ subdividing $\s_0$, then
there exists a unique minimal 
regular subdivision $\F^{\reg} $ of $\F$, that is,  if $\F'$ is any other regular subdivision 
of $\F$ then $\F'$ subdivides  $\F^{\reg} $ (see \cite[Section 10.2]{Cox}).  If 
$\theta' = \R_{\geq 0} \, u+ \R_{\geq 0}  \, v  \in  \F^{\reg}$ is a two dimensional cone, 
with  
$u= a_1 e_1 + a_2 e_2$ and $v = b_1 e_1 + b_2 e_2$, 
then  $a_1 b_2 - a_2 b_1 = \pm 1 $ and
the chart 
$\psi^{\theta}_{\sigma_0}: X_\theta \to X_{\s_0}$ 
of the modification 
\begin{equation}  \label{Pi-sigma}
\psi^{\F^\reg}_{\s_0} : X_{\F^\reg}  \to X_{\s_0} 
\end{equation}
is the monomial map given by:
\begin{equation} \label{pi-sigma}
\begin{array}{lcl}
x &=x_\theta ^{a_{1}} y_\theta ^{b_{1}}, 
\\
y &=x_\theta^{a_{2}} y_\theta^{b_{2}}.
\end{array}
\end{equation}
If $\check{u}, \check{v} \in M$ denotes the dual basis of 
$u, v$ then,  one has 
$x_\theta := \chi^{\check{u}}, y_\theta := \chi^{\check{v}}$ and
$\C[X_\theta] = \C[ x_\theta, y_\theta]$. 
Notice that the closure of the orbit $O_\r$ associated with 
a ray $\rho = \R_{\geq 0} u$, 
is defined on the chart \eqref{pi-sigma} by $ x_\theta =0$.

\subsection{Newton polygons, Newton fans and support functions} \label{toric}

Let $f = \sum a_{i, j} x^i y ^j  \in \C[[x, y]]$ be a nonzero power 
series. The \textbf{support}  $\boxed{\mathrm{supp} (f)}$ of the power series $f$ consist of those vectors 
$(i, j) = i \check{e}_1 + j \check{e}_2 \in  M$ with nonzero coefficient $a_{i, j}$. 
Recall that   $\check{\sigma}_0 =  \R_{\geq 0} \check{e}_1 + 
\R_{\geq 0} \check{e}_2$ is the dual cone of $\sigma_0 = \R_{\geq 0} {e}_1 + 
\R_{\geq 0} {e}_2 $ (see the  notation of Section \ref{sec:fan}). 

\medskip 

The \textbf{Newton polygon} $\boxed{\mathcal{N} (f)}$ is the convex hull of $\mathrm{supp} (f) + \check{\s_0}$. 
The face $\mathcal{E}_u$ of  $\mathcal{N} (f)$ defined by a vector $u \in \s_0$ is the set 
of  elements $v \in \mathcal{N} (f)$ such that 
$
\langle v,u \rangle = 
	\min \{  \langle v',u \rangle \mid   v' \in \mathcal{N} (f)   \}$.
		The face  $\mathcal{E}_u$ is compact when $u$ belongs to 
	the interior of $\s_0$. All the faces of $\mathcal{N} (f)$ are of this form.
If $\mathcal{E}$ is a face of $\mathcal{N} (f) $ the closure  of the set 
$\{ u \in \s_0 \mid \mathcal{E}_u = \mathcal{E} \}$ is a cone $\theta_{\mathcal{E}} \subset \sigma_0$ which is 
rational for the lattice $N$. 
The
set  $\boxed{\mathcal{F} (f)} $ consisting  of cones $\theta_{\mathcal{E}}$, for $\mathcal{E}$ running through the 
faces of $\mathcal{N}(f)$, is a fan of $N$ supported on $\s_0$, called the \textbf{Newton fan} of $f$. 
We denote by $\boxed{\F^\reg(f)}$ the minimal regular subdivision of the fan $\mathcal{F} (f) $.

\medskip 
The \textbf{support function} $\boxed{\Phi_{\mathcal{N}}} : \sigma_0 \longrightarrow \R_{\geq 0}$ 
of the polyhedron $\mathcal{N} := \mathcal{N} (f) $  is defined by
$
\Phi_{\mathcal{N}} (v) =\underset{u \in \mathcal{N} } {\mathrm{min}} \langle v, u \rangle$.
The support function $\Phi_{\mathcal{N}}$ is linear precisely on each cone of the fan $\F (f) $. 
If $\theta \in \F (f) $ is a two dimensional cone 
then there exists a unique vertex $u_\theta$  of the Newton polygon $\mathcal{N}(f) $
such that 
$\Phi_{\mathcal{N}} (v)  = \langle v, u_\theta \rangle $ for all $v \in \theta$.
In addition, for any real number $\xi > 0$ we have that 
\begin{equation} \label{homothetic}
\Phi_{\xi \mathcal{N}} =  \xi  \Phi_{\mathcal{N}}.
\end{equation}
Notice that the notion of support function can be defined for any convex polyhedra and 
determines it  (see \cite[Theorem 3.8, and Theorem 6.8]{EwCombConv}). In particular, 
we have 
\begin{equation} \label{eq: int-semi}
\mathcal{N} (f) 
= \cap_{v \in \F(f) _\prim} \{ u \in M_\R \mid \langle v , u \rangle  \geq \Phi_{\mathcal{N}} (v) \} =  \cap_{v \in \sigma_0} \{ u \in M_\R \mid \langle v , u \rangle  \geq \Phi_{\mathcal{N}} (v) \}.
\end{equation}
\medskip

A vector $v \in \sigma_0$  defines a \textbf{monomial valuation} $\boxed{\ord_v}$ of the completion  
 $\hat{\mathcal{O}} \cong \C[[ x, y ]]$  of the local ring $\mathcal{O}$,
given by 
\begin{equation} \label{monomial_val}
\ord_{v} (h) : = \min \{ \langle v, u \rangle \mid v \in \mathrm{supp} (h) \}, \mbox{ for }  h \in  \C[[ x, y ]] \setminus \{0\}.
\end{equation}
It follows that 
for any $v \in \sigma_0$  and any plane curve $C$
we have 
\begin{equation} \label{eq: sf}
\ord_v (C) = \Phi_{ \mathcal{N} (f_C) }  (v). 
\end{equation}

 The following  lemma  will be useful later.
\begin{lema} \label{lem :inclusion}
Let $A$ be a non-empty subset of $\check{\s}_0$, set 
$\mathcal{A} := A + \check{\s}_0$, consider  a rational number $\xi >0$, and 
let $\mathcal{G}$ be subdivision of $\F (f)$. 
The following conditions are equivalent. 
\begin{enumerate}[label=({\alph*})]
\item \label{inc-a}
$\mathcal{A}$ is contained in the interior of the polygon $\xi \mathcal{N}(f) $,
\item \label{inc-b}
$\Phi_{\mathcal{A}} (v) >  \xi \Phi_{\mathcal{N}(f)} (v) $ 
for any $ v \in \mathcal{G}_\prim$,
\item \label{inc-c}
$\Phi_{\mathcal{A}} (v) >  \xi \Phi_{\mathcal{N}(f)} (v) $,
for any $ v \in \F (f) _\prim$.   
\end{enumerate}
\end{lema}
\begin{proof}
By \eqref{homothetic} and  \eqref{eq: int-semi}, 
the interior of $\xi \mathcal{N}(f) $ is the intersection of open half-spaces
\[
\bigcap_{v \in \F(f)_\prim} \{ u \in M_\R \mid \langle v , u \rangle  >  \xi \Phi_{\mathcal{N}(f)} (v) \} =  \bigcap_{v \in \sigma_0} \{ u \in M_\R \mid \langle v , u \rangle  >  \xi \Phi_{\mathcal{N}(f)} (v) \}.
\]
This implies the equivalence between the conditions \ref{inc-a}, \ref{inc-b} and \ref{inc-c}. 
\end{proof}

\subsection{Newton modifications} 
We present the properties of Newton modifications 
associated with a germ of plane curve on a smooth surface relative to a cross
following \cite{GBGPPP20}.

\medskip

Let $(R, L)$ define a cross \index{cross} at a point $o$ of a smooth surface $S$.  
The set of divisors supported on $R + L$ is a rank two lattice $\boxed{M_{R, L}}$
 with basis $R, L$. The map $M_{R, L} \to M$, which sends $aR + b L $ to $a  \check{e}_1 + b 
\check{e}_2$  is 
an isomorphism of lattices, and 
it extends to an isomorphism of 
real vector spaces which maps the cone  $\check{\s}_0^{R,L}$ of real effective divisors supported on $R + L$
onto the cone $\check{\s}_0$. 
We denote by $\boxed{N_{R, L}} \simeq N$ the dual lattice of $M_{R, L}$, by $\boxed{e_{R}, e_{L}} \in N_{R, L}$ the dual basis of $R =: \check{e}_{R}, L =: \check{e}_{L}$, 
and by $\boxed{\s_0^{R,L}}$ the cone $\R_{\geq 0} e_{R} + \R_{\geq 0} e_{L}$.

\medskip

Let us consider the projectivization map: 
\[
\begin{array}{ccc}
\boxed{\phi}: (N_{R,L})_\R \setminus \{ 0 \} & \longrightarrow & \mathbb{P}\left((N_{R,L})_\R\right)
\\
a  e_R + b e_L & \mapsto & (a: b), 
\end{array}
\]
where $(a:b)$ denote homogeneous coordinates.
If $u = a e_R + b e_L $ is nonzero the \textbf{slope}
$\boxed{\slp (u)}: = b/a \in {\R} \cup \{  \infty \}$
is the affine coordinate of the point $\phi (u) \in \mathbb{P}\left((N_{R,L})_\R\right) \cong \mathbb{P}^1_\R$.

\begin{notation} \label{not:proj}
If $a, b \in \N$ are coprime we abuse of notation by denoting 
 the image of the vector $u = a e_R + b e_L $ by $\phi$ also by $u$.
 \end{notation}

\begin{defi} \label{slopetrunk}
 Let $\F$ be a fan  of $N_{R, L}$ subdividing $\s_0^{R, L}$. 
 The  {\bf trunk} $\boxed{\theta(\fan)}$ of $\F$  is 
 the segment $\phi ( \s_0^{R, L} \setminus \{ 0 \} ) =[e_R, e_L] $ with the finite set of marked points defined by the image by $\phi$ of the rays of $\F$. It is endowed with the slope coordinate function
 \[
 \slp: [e_R, e_L] \to [0, \infty] \subset \mathbb{P}^1_\R.
 \]
\end{defi}
For example, with Notation  \ref{not:proj} we have that if $\F$ is the subdivision of $\s_0^{R, L}$ with rays spanned by 
the vectors $e_R$, $e_R + e_L$ and $e_L$ the trunk $\theta({\F})$ is 
the segment $[e_R, e_{L}]$ with marked points $e_R$, $e_R + e_L$ and $e_L$
which have slope coordinate $0$, $1$ and $\infty$ respectively.

\medskip 

Any  regular fan $\F$ with respect to $N_{R, L}$ subdividing $\s_0^{R,L}$
defines a model of $(S,o)$:
\begin{equation} \label{toroidal-map}
\boxed{\psi^\F_{R, L}}:  S_{\F} \to S, 
\end{equation}
which is 
defined by gluing monomial maps of the form $\psi^{\theta}_{\sigma_0}$ 
for $\theta \in\F$.

\medskip

\begin{remark} \label{rem:bj}
There is a natural bijection between the set of marked points of the trunk
$\theta(\fan)$ and the set of irreducible components of 
$(\psi^\F_{R, L})^{-1} (R +L)$ which sends 
a marked point $p$ to the irreducible component $D_p$
    of $(\psi^\F_{R, L})^{-1} (R +L)$ which 
    contains the orbit 
    labelled by the ray of $\F$ of slope $\slp (p)$. 
    In particular, if  $p = e_R$ (resp. $p= e_L$) 
then one has that $D_p = R$ (resp. $D_p =L$). 
\end{remark}
\begin{notation} \label{not:mp}
Using this bijection explained above, 
the marked points of $\theta(\fan)$
can be relabelled as 
$e_{D}$, for $D$ running through the irreducible components of 
$(\psi^\F_{R, L})^{-1} (R +L)$. 
\end{notation}

\begin{remark} \label{van-val} 
If $v \in  \s_0^{R,L}$, then the monomial valuation $\ord_v$ defined in 
\eqref{monomial_val} is independent of the choice of local coordinates $(x, y)$ defining 
the cross $(R,L)$. In the particular case of the vectors $e_R$ and $e_L$,  we get that
$\ord_{e_R}  = \nu_{R}$,   and $\ord_{e_L} = \nu_{L}$
are the vanishing order valuations  along the branches $R$ and $L$ respectively.
If $v = n e_R + m e_L \in N_{R, L}$ for $n, m \in \N$ coprime,  then  
the monomial valuation  $\ord_{v}$ is a divisorial valuation. In order to see this, take any regular fan $\F$ of $N_{R, L}$ subdividing $\s_0$ and containing the ray $\R_{\geq 0} v$. 
Then, if $p \in (e_R, e_L)$ has slope $m/n$ then we have that 
$\ord_v = \nu_{D_p}$.
\end{remark}

Let $(x, y)$  be a pair of local coordinates defining the cross  $(R, L)$. 
Let  $C$ be a plane curve singularity defined by a power series $f  \in \C[[x, y]]$.
The \textbf{Newton polygon of $C$ with respect to the cross} $(R, L)$  
is 
$\boxed{\mathcal{N}_{R, L} (C)} $ is just the polygon $\mathcal{N} (f)$ seen as a subset of $(M_{R, L})_\R$.
The \textbf{Newton fan of $C$ with respect to the cross} $(R, L)$ is denoted by
$\boxed{\F_{R, L} (C)}$ is just the fan $ \F (f)$, whose support  is seen as a subset of $(N_{R, L})_\R$. 
The polygon $\mathcal{N}_{R, L} (C)$ and the fan $\F_{R, L} (C)$ 
are independent of the choices of local coordinates $(x, y)$ defining 
the cross $(R,L)$ and the function $f_C \in \C[[x, y ]]$ defining $C$ (see \cite[Sec. 4.1]{GBGPPP20}).

\begin{defi}  \cite[Def. 4.14]{GBGPPP20}
We denote by $\boxed{\F^\reg_{R, L} (C)}$ the minimal regular subdivision of the fan $\F_{R, L} (C)$.
The map 
\begin{equation} \label{rnm}
\boxed{\psi_{R, L}^{C, \reg}}:= \psi^{\F^\reg_{R, L} (C)}: S_{\F^\reg_{R, L} (C)} \to S 
\end{equation}
is called the \textbf{regularized Newton modification} of $S$  defined by $C$ with respect to the cross $(R,L)$.
\end{defi}

Notice that the  map defined by \eqref{rnm} is a model of $(S,o)$, 
in particular, its exceptional divisor has simple normal crossings.
Some concrete cases are discussed in Example \ref{ex:2branches} below.

\begin{remark}
In \cite{GBGPPP20}, a more general notion of   \textit{Newton modification} 
$
{\psi_{R, L}^{C}} := \psi^{\F_{R, L} (C)}: S_{\F_{R, L} (C)} \to S$, is considered. Notice that the surface $S_{\F_{R, L} (C)}$ may be singular.
\end{remark}

\subsection{Toroidal resolutions of plane curves}\label{subsec:torres}

In this section we summarize an algorithm of toroidal embedded resolution given in  \cite[Section 4]{GBGPPP20}.

\medskip 

\begin{defi} \cite[Def. 3.29 and 4.15]{GBGPPP20}  \label{def:threeres}
A (smooth) toroidal surface is a smooth complex analytic surface $\Sigma$ endowed with a
 normal crossing divisor $\boxed{\partial \Sigma}$, called its \textbf{boundary}.  A modification $\pi: (\Sigma_2, \partial \Sigma_2) \to (\Sigma_1, \partial \Sigma_1)$ between toroidal surfaces is a
\textbf{toroidal modification} if $\pi^{-1} (\partial \Sigma_1) \subset \partial \Sigma_2$.

  \medskip 
  
Let $C$ be a plane curve singularity 
on a germ of smooth surface $(S, o)$ 
endowed with a germ of normal 
crossing divisor $\partial S$, and 
let $(\Sigma, \partial \Sigma)$ be a smooth toroidal surface. 
A toroidal modification $\pi: \Sigma \to S$ is a 
\textbf{toroidal embedded resolution  of $C$} if
the boundary $\partial \Sigma$ of $\Sigma$ contains the 
     reduction of the total transform $\pi^*(C)$ of $C$ by $\pi$.
 The reduction of the image $\pi(\partial \Sigma)$ of $\partial \Sigma$ in $S$
    is called the {\bf completion $\boxed{\hat{C}_{\pi}}$ of $C$ relative to $\pi$}.     
  \end{defi}

 If $\pi$   is a toroidal embedded resolution of $C$, 
 then the strict transform of $C$ by $\pi$ is smooth and transversal to the exceptional divisor of $\pi$, 
 hence it is an embedded resolution of $C$. 
  The simplest case of toroidal embedded resolution
 is given by the following example. 
     \begin{example}
     If $\pi: S_\pi \to S$ is the minimal embedded resolution of $C$ then 
     $\pi$ is a {toroidal embedded resolution  of $C$},  where
     $\partial S := \emptyset$ and $\partial S_\pi$  is the reduced divisor 
    of $\pi^{*} (C)$. In this case the completion $\hat{C}_\pi$ is the reduction of $C$. 
     \end{example}
     
      We focus from now on toroidal modifications  
$\pi: \Sigma \to S$, where 
$(R, L)$ is a fixed cross on $(S,o)$ and  $\partial S := R + L $.
 
\begin{example}
Denote by $\psi:  S_\psi \to S$ 
the regularized Newton modification \eqref{rnm} of $S$ defined by 
$C$ with respect to the cross $(R,L)$. 
Then, $\psi$ is a modification of toroidal surfaces,  when we take $\partial S := R+ L$ and  $\partial S_\psi$ is any reduced normal crossing divisor on $S_\psi$
containing the reduction of $\psi^*(R+ L)$.
\end{example}

\begin{prop}{\cite[Prop. 4.18]{GBGPPP20}}\label{ToricModif} 
Let
 $(C,o)$ be a plane curve singularity 
  on a smooth surface $S$, and let $(R, L)$ be a cross at $o$.
Assume that neither  $R$ nor $L$ is a branch of $C$. 
Denote by $\psi:  S_\psi \to S$ 
the regularized Newton modification \eqref{rnm} of $S$ defined by 
$C$ with respect to the cross $(R,L)$. 
Then, the strict transform $C^\psi$ of $C$ 
intersects the reduced divisor of $\psi^* (R +L)$
only at smooth points of it.
\end{prop}

By  Proposition \ref{ToricModif}, one of the following two cases holds at 
 each point of intersection $o_i$ of the strict transform $C^\psi$ with the exceptional divisor 
of $\psi$. 
   \begin{enumerate}
\item \label{st1} The germ of the strict transform $C^\psi$ at the point $o_i$ is smooth and 
transversal to the exceptional divisor of $\psi$.  Then,  only one branch of $C^\psi$ passes through $o_i$,  and together with 
the germ $R_i$  of the exceptional divisor define a canonical cross on $S_\psi$.

\item \label{st2} Otherwise,  we can choose 
a smooth germ $L_i$ transversal to 
the germ $R_i$ of exceptional divisor at $o_i$ and then
$(R_i, L_i)$ defines a cross at $o_i$.
\end{enumerate}
In case \eqref{st2} we have a plane curve singularity, the germ of strict transform $C^\psi$ of $C$ at $o_i$, 
and a cross at $o_i$, so that we can apply 
to it the associated 
regularized Newton modification defined by it with respect to this cross. 
 This leads to the following 
algorithm.

\begin{algorithm}{\cite[Alg. 4.22 and Prop. 5.1]{GBGPPP20}}  \label{alg:tores}  
    Let $(S,o)$ be a smooth germ of surface, $R$ be a smooth branch on $(S, o)$,   
    and $C$ be a reduced germ of curve on $(S,o)$, which does not contain the 
    branch $R$ in its support. 
  \medskip
  
  \noindent
   {\bf STEP 1.}  If $(R, C)$ is a cross, then STOP. 
      
     \noindent
    {\bf STEP 2.}  Choose a smooth branch $L$ on $(S,o)$, possibly 
             included in $C$, such that $(R, L)$ is a cross.  

   \noindent 
  {\bf STEP 3.}  
  Consider 
the regularized Newton modification $\psi:= \psi_{R, L}^{C, reg}: (S_\psi ,  \partial S_\psi) \to (S,\partial S)$
                 of $S$ defined by $C$
                 with respect to the cross $(R, L)$, 
                 where $\partial S := R+ L$ and $\partial S_\psi:= \psi^{-1} (R+ L)$, 
                 and the strict transform $C^\psi$ 
                 of $C$ by $\psi$. 
 \medskip

    \noindent
   {\bf STEP 4.} For each point $\tilde{o}$ belonging to 
            $C^\psi \cap  \partial S_\psi$, denote:
                     \begin{itemize}
                          \item $R:=$ the germ of $\partial S_\psi$
                         at $\tilde{o}$; 
                          \item $C:=$ the germ of $C^\psi$ at $\tilde{o}$; 
                          \item $ o := \tilde{o}$;
                          \item  $S := S_\psi$. 
                     \end{itemize}

        \noindent
      {\bf STEP 5.}  Go to step 1.              
\end{algorithm}

\begin{prop}{\cite[Prop. 5.1]{GBGPPP20}}  \label{Prop: toroidal res}
The  Algorithm \ref{alg:tores} stops after finitely many iterations and 
provides a toroidal embedded resolution $
\pi: (\Sigma, \partial \Sigma) \to (S, \partial S) 
$ of $C$ with the following boundaries. 
\begin{itemize}
\item
 $\partial S := R + L$,  where $L$ is the branch fixed at the first the first iteration of 
STEP 1 or STEP 2 of the Algorithm \ref{alg:tores}. 
\item 
 $\partial \Sigma$ is the reduced normal crossing divisor 
which contains the reduction of $\pi^* (C)$ and 
the strict transforms of the components of  
the crosses considered when running the Algorithm \ref{alg:tores}.
\end{itemize}

\end{prop}

\begin{example} \label{ex:2branches}
Let $(x, y)$ be a local coordinate system on the surface $(S, o)$ 
and consider $f_1=(y^2+x^3)^2+x^6y, f_2=y^3+x^5 \in \mathcal{O}$.
We describe a toroidal embedded resolution of the plane curve singularity
$C = C_{f_1} + C_{f_2}$ following Algorithm \ref{alg:tores}. First, we fix the smooth branch $R_1 = C_x$. 
At the second step of the algorithm we choose $L_1 = C_y$, 
and the cross $(R_1, L_1)$ at $o_1$. 
The Newton polygon $\mathcal{N}_{R_1,L_1}(C)$ has two compact edges 
which are orthogonal to the vectors $2 e_{R_1} + 3 e_{L_1}$ and $3 e_{R_1} + 5 e_{L_1}$. 
These vectors span rays of the fan $\F_1 :=  \F^\reg_{R_1, L_1} (C)$, which is represented
in the left part of  Figure \ref{fig:NewtonFan1}. 

\begin{figure}[!ht]
\centering
\includegraphics[width=.8\linewidth]{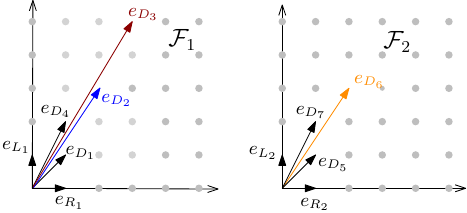}
\caption{The regularized Newton fans $\F_1$ and $\F_2$.}
\label{fig:NewtonFan1}
\end{figure}

Let us take the regularized Newton modification 
 $\psi_1 := \psi_{R_1, L_1}^{C, \reg}: S_{\F_1} \to S$.
By Proposition \ref{ToricModif}, the strict transform 
$C^{\psi_1}$ of $C$ intersects the exceptional divisor 
of $\psi_1$ at  points $o_2 \in D_2$ and $o_3 \in D_3$, where the labels are those of Figure \ref{fig:NewtonFan1}. 
We can check this on  the  chart 
$x = x_2^2 \, x_3^3$,   $y =x_2^3 \, x_3^5$ of $\psi_1$, where 
$D_2 = C_{x_2}$ and $D_3 = C_{x_3}$, and 
the total transform of $C_1$ (resp. of $C_2$) is 
defined by $x_2^{12} \, x_3^{18} \, \left((x_3+1)^2 + x_2^3  \, x_3^5 \right) =0$
(resp. by $x_2^9  \, x_3^{15}(1+x_2) =0$). 
Thus, the point  $o_2$ (resp. $o_3$) has coordinates $(x_2, x_3) = (0, -1)$ 
(resp. $(x_2, x_3) = (-1, 0)$).  Then, we iterate the algorithm at the points $o_3$ and $o_2$: 

- At the point $o_3$ we 
get the cross 
$(R_3: = D_3 , L_3 := C^{\psi_1} = C_{2}^{\psi_1})$,
and we stop at step 1. 

- At the point $o_2$, we choose  the cross $(R_2 := D_2, L_2 := C_{x_3 +1} )$ at the second 
step of the algorithm. 
 The fan $\F_2 :=  \F^\reg_{R_2, L_2} (C)$, is represented
in Figure \ref{fig:NewtonFan1}.
It defines 
the regularized Newton modification 
 $\psi_2 := \psi_{R_2, L_2}^{C, \reg}: S_{\F_2} \to S_{\F_1}  $.
By Proposition \ref{ToricModif}, 
the strict transform of $C$ intersects the 
exceptional divisor of $\psi_2$ at 
a point $o_4 \in D_6$ (the labels of the components of the exceptional divisor 
are indicated in the right part of Figure \ref{fig:NewtonFan1}).
We obtain that 
$\pi := \psi_1 \circ \psi_2$
is a toroidal embedded resolution of $C$, since 
at the point $o_4$ we get the cross 
$(R_4 = D_6, L_4 := C^\pi = C_1^\pi) $,
that is,  the iteration of the algorithm at the point $o_4$ stops at step 1. 
The image of $L_2$ on the inicial surface is the branch 
$C_{y^2 + x^3}$. 
\end{example}

\subsection{The fan tree of the toroidal embedded resolutions of Algorithm \ref{alg:tores}}

We explain now how a tree, called the \textit{fan tree},  can be associated with a toroidal embedded 
resolution of the form given in Proposition \ref{Prop: toroidal res}, following 
\cite{GBGPPP20}.
The fan tree encodes the combinatorial structure of 
the toroidal resolution process. The fan tree is a variant of 
the \textit{Newton tree} considered by Cassou-Nogu\`es and Libgober
with equivalent decorations  (see  \cite{CNL14}). 
We introduce first some notations.

\begin{notation}   \label{def:manycrosses}
Assume that one executes Algorithm \ref{alg:tores}  on $(S,o)$, starting from the curve singularity $C$ and the smooth branch $R$, 
which is not a component of $C$. 
We denote by $\boxed{I_\pi}$,  a finite set labelling the infinitely near points $o_i$ of $o$ at which one applies STEP 1 or STEP 2. We  assume that $1 \in I_\pi$ and then $o_1 =o$ and $R_1 = R$. If $i \in I_\pi $ we denote by $(R_i, L_i)$ the corresponding cross at $o_i$. 
Denote by $\boxed{ e_{R_i}, e_{L_i} }$ the canonical basis of the 
weight lattice $N_{R_i, L_i}$.  If $i \in I_\pi$ and 
$i \ne 1$, the branch $R_i$ is included in the exceptional divisor 
of the regularized Newton modification performed at the previous step.  
We denote by $L_{i, \pi}$, or simply by $L_i$, the projection on $S$ of the curvetta $L_i$ at $o_i$. 
We denote by $\psi_i := \psi_{R_i, L_i}^{C, reg}$ the regularized Newton modification of $C$ with respect to the cross $(R_i, L_i)$ at the point $o_i$.  We consider the trunk  
$\theta_i$ of the fan $\fan_{R_i, L_i}^{\reg}(C)$, for  $i \in I_\pi$. The trunk $\theta_i$ is 
the segment  $[e_{R_i}, e_{L_i}]$ endowed with 
the slope coordinate function  $\slp_i : 
[e_{R_i}, e_{L_i}] \to [0, + \infty]$, 
and with marked points defined by the edges of the fan $\fan_{R_i, L_i}^{\reg}(C)$ (see Definition \ref{slopetrunk}). 
Recall that we label the marked points  of the trunk $\theta_i$ by $e_D$, where $D$ runs through the irreducible components of the reduction of the divisor $\psi_i^* (R_i + L_i)$ (see Notation \ref{not:mp}).
By definition, the strict transform of $D$ on the final surface $\Sigma$, which we denote also by $D$, is a component of the boundary $\partial \Sigma$. Notice that $D$ is an exceptional prime 
divisor of $\pi$ precisely when  $e_D \in (e_{R_i}, e_{L_i})$, or 
$D = R_i$ for $i \ne 1$. Otherwise $D \subset \partial \Sigma$ is the strict transform of a component of $\hat{C}_\pi = R + \sum_{j \in I_\pi} L_{j,\pi}$. 
\end{notation}

\begin{remark} 
If $C_j$ is a component of $C$ then there exists a unique index $i_j \in I_\pi$ such that $C_j = L_{i_j,\pi}$. 
 This means that the algorithm stops at the point $o_{i_j}$, since it is the point of intersection of the strict transform of $C_j$ with the exceptional divisor $E_\pi$. In this case 
the morphim $\psi_{i_j}$ is the identity map and 
the marked points of the fan tree $\theta_{i_j}$ are $e_{R_{i_j}}$ and $e_{C_j}$.                     
\end{remark}

\begin{defi}  \label{def:fantreetr}      The {\bf fan tree $\boxed{\theta_{\pi}(C)}$ of the toroidal embedded resolution  
    $\pi: \Sigma \to S$ of $C$} of Proposition \ref{Prop: toroidal res}
    is a tree 
 endowed with a finite set of marked points.
The set   $\theta_{\pi}(C)$ is obtained from 
        the disjoint union of the trunks
        $\theta_i = \theta(\fan_{R_i, L_i}^{\reg}(C))$, for ${i \in I_\pi}$, by identifying
     the marked points  labeled by the same 
        irreducible component of $\partial \Sigma$. 
\end{defi}
We denote in the same way each interval $[e_{R_i}, e_{L_i}]$ and its image in $\theta_{\pi}(C)$.
By definition, if $e \in \theta_{\pi}(C)$ and $e   \ne e_R$ there exists a unique index $i \in I_\pi$ such that $e \in (e_{R_i}, e_{L_i}] $.
This property allows to endow the fan tree with 
its slope function:
$
\boxed{\slp_{\pi}} :  \theta_{\pi}(C) \to [0, \infty]
$,
defined by 
            \[
             \slp_{\pi}(e) = \left\{ 
             \begin{array}{ccl}
             0 & \mbox{ if } & e = e_R,
             \\
             \slp_{i} (e)  & \mbox{ if } & e \in
             (e_{R_i}, e_{L_i}] \mbox{ for some } i \in I_\pi.
             \end{array}
             \right.
             \]
This function is not continuous precisely on the set
ramification points 
$\{e_{R_i} \mid i \in I_\pi, i \ne 1 \}$  
of the 
tree $\theta_\pi (C)$.

 \medskip

By construction, the toroidal embedded resolution $\pi$ is also an embedded resolution 
of the completion $\hat{C}_\pi$, and the boundary $\partial \Sigma$ is 
the reduced divisor of the total transform of $\hat{C}_\pi$. 
The associated dual graph is determined in terms of the fan tree.

\begin{prop} {\cite[Prop. 4.35]{GBGPPP20}} 
 \label{prop:fantreedualgr}
   The fan tree $\theta_{\pi}(C)$ is isomorphic to the dual graph
   $G (\partial \Sigma)$ of 
   the boundary $\partial \Sigma$ of the source of the toroidal embedded resolution 
   $\pi:\Sigma \to S$ of $C$
   of Proposition \ref{Prop: toroidal res}  
   by an isomorphism which respects the labels. 
\end{prop}

\begin{remark} $\,$ \label{rem:tree}
 By Proposition  \ref{prop:fantreedualgr}
marked points of $\theta_\pi (C)$ which are of valency $\geq 2$ 
are labelled by the
   irreducible components of the exceptional divisor of $\pi$.
   The ramification points of $\theta_\pi (C)$ 
   are labeled by the elements of $\mathcal{R}_\pi (C)$, 
   that is, by the 
   rupture components of the divisor $(\pi^* (\hat{C}_\pi))^{red} = \partial \Sigma$. 
   The end points of the tree $\theta_{\pi}(C)$, which correspond 
   to the vertices of $G (\partial \Sigma)$ of valency one, 
   are labeled by 
   the irreducible components of the completion $\hat{C}_\pi$.
   This implies that the set of components of 
   $\hat{C}_\pi$ contains a set of maximal contact curves of the pair $(\pi, C)$, see Definition \ref{mcc}.
\end{remark}

\begin{remark}  \label{rem:minimal}
The minimal embedded resolution of $C$ can be obtained as a toroidal embedded resolution by choosing a suitable reference smooth branch $R$ and suitable auxiliary branches $L_i$ at the second step of the Algorithm \ref{alg:tores}.
One may take a \textit{maximal contact toroidal embedded resolution}, see  \cite[Sec. 4.2]{PhD-Robredo}. It is also the case of some
toroidal resolutions described in \cite{LO}.
\end{remark}

\begin{example} \label{Ex:2BranchesFanTree}
We describe the fan tree of the toroidal embedded resolution of Example \ref{ex:2branches}. 
In the left side of Figure \ref{fig:2BranchesFanTree} we have represented 
the trunks $\theta(\fan_{R_i, L_i}(C))$ for $i = 1,\dots, 4$. 
\begin{figure}[!ht]
\centering
\includegraphics[scale=0.9]{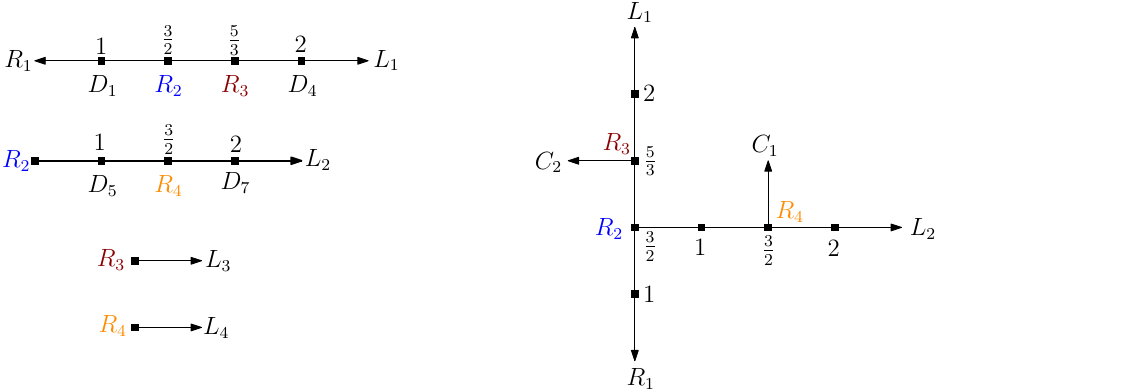}
\caption{The trunks and the  fan tree $\theta_\pi(C)$ of Example \ref{ex:2branches}.}
\label{fig:2BranchesFanTree}
\end{figure}
The fan tree $\theta_{\pi} (C)$, which is represented in right part of Figure \ref{fig:2BranchesFanTree}, 
is obtained from these trunks by glueing the points with the same label. 
The completion 
of the toroidal embedded resolution $\pi$ is 
 $\hat{C}_\pi = R_1 + L_1 + L_2 + C_1 + C_2$.  
One may check that this is the minimal embedded resolution of $C$.
\end{example}

The following notion of representing divisor is equivalent to the notion of
\textit{representing divisor of a rational point of an Eggers-Wall tree} in \cite{GGPValTree}, 
see  \cite[Section 1.6]{GBGPPP20}.

\begin{defi} \label{def:repdiv}
A point $p \in \theta_{\pi}(C)$ is \textbf{rational} if $ \slp_{\pi} (p) \in \Q^*$.
The \textbf{representing divisor} $\boxed{D_p}$ of a rational point $p \in \theta_{\pi}(C)$  is defined as follows.
By Definition \ref{def:fantreetr}  
there exists a unique $i \in I_\pi$ such that $p \in (e_{R_i}, e_{L_i})$ and then $\slp_{i} (p)  = \frac{m}{n}$, 
where $n, m >0$ are two coprime integers. 
Take a regular subdivision  $\F$ of the fan $\fan_{R_i, L_i}^{\reg}(C)$ which contains the ray spanned by the vector 
\begin{equation}\label{eq:eP}
\boxed{ e_{D_p}} :=  n e_{R_i} +  m e_{L_i} \in N_{R_i, L_i}.
\end{equation}
Then, the prime exceptional divisor $D_p$ is an irreducible component of 
on the source of the toric model 
$\psi^{\F}_{R_i, L_i}$ (see Remark \ref{rem:bj}).
\end{defi}  

We can replace the model $\psi_i$ by the map $\psi^{\F}_{R_i, L_i}$ considered in Definition \ref{def:repdiv}, in the running of  Algorithm \ref{alg:tores}.
Then, the output is a model $\pi': (\Sigma', \partial \Sigma') \to (S, \partial S)$ of $C$ dominating $\pi: (\Sigma, \partial \Sigma) \to (S, \partial S)$,
such that  the exceptional divisor  
$D_p$ appears as a component of $\partial \Sigma'$. 
In particular, the representing 
divisor $D_p$ of a rational point $p$ appears on $\Sigma$
if and only if $p$ is a marked point $p$  of valency $\geq 2$ of $ \theta_\pi (C)$.

\section{Multiplier ideals and Newton polygons} \label{sec:main}

In this section we fix a plane curve singularity $C$, a smooth branch $R$  
which is not a component of $C$ and 
a toroidal embedded resolution $\pi$ of $C$ given by Proposition 
\ref{Prop: toroidal res}. 
Recall that the set $I_\pi$ indexing  the crosses appearing in the Algorithm \ref{alg:tores}
was introduced in Notation \ref{def:manycrosses}.

\begin{defi}
  The \textbf{log-discrepancy vector at $o_i$}, for $i\in I_\pi $ is 
    \[
    \boxed{ \underline{{\l}}_i } :=  \l_{R_i} R_i  +    L_i \in M_{R_i, L_i}.
     \] 
\end{defi}
       Let $D$ be a prime component of $\hat{C}_{\pi}$. Recall that if $D$ is exceptional then
      $\lambda_{D}$ denotes the log-discrepancy $D$. If $D$ is not exceptional we set $\l_{D} :=1$.

\begin{notation}
Let $A$ be a curve on $S$. If $i \in I_\pi$ we 
denote by  $\boxed{\Newton_{R_i, L_i }(A)}$
                    the Newton polygon  of the germ  of  the total transform of $A$ at $o_i$
                  relative to the cross $(R_i, L_i)$, 
                    and by 
                     $\boxed{\fan_{R_i,  L_i }(A)}$ the
                     corresponding Newton fan.
                     Recall that
                     $\Phi_{\mathcal{N}_{R_i, L_i} (A) }$
                     denotes the support function of 
                     the Newton  polygon $\mathcal{N}_{R_i, L_i} (A) $
                     (see Section \ref{toric}).

\end{notation}

In the following proposition we apply the notion of representing divisor $D_p$ of a rational point $p$ in the fan tree (see Definition 
\ref{def:repdiv}). 
\begin{prop}  \label{mcp} 
Let $p$ be a rational point of the fan tree $\theta_\pi (C)$,  
let $i \in I_\pi$ be the unique index such that $p \in (R_i, L_i)$,
and consider the vector  $e_{D_p}  \in N_{R_i, L_i}$ defined by \eqref{eq:eP}.  
Then, for any plane curve $A$ on $S$ 
we have  
\begin{equation} \label{eq: sf2}
\nu_{D_p} (A) = \Phi_{\mathcal{N}_{R_i, L_i} (A) } ( e_{D_p}),  
\end{equation}
\begin{equation} \label{eq: sf4}
\nu_{L_{i, \pi}} (A) = \Phi_{\mathcal{N}_{R_i, L_i} (A) } (e_{L_i}).   
\end{equation}
In addition,  the log-discrepancy of the exceptional prime $D_p$ is given by
\begin{equation} \label{eq: sf5}
\lambda_{D_p} = \langle  e_{D_p} , \underline{{\l}}_i \rangle. 
\end{equation}
\end{prop}
\begin{proof}
Take local coordinates $(x_i, y_i)$ defining the cross $(R_i, L_i)$ at $o_i$. 
Since $p$ is rational, we have that $\slp_{\pi} (p)  = \frac{m}{n}$, 
where $n, m >0$ are coprime. The
monomial valuation 
$\ord_{ e_{D_p}}$ is defined on the completion of the local ring at $o_i$ 
(identified with $\C [[ x_i, y_i]]$),  by 
\[
\ord_{ e_{D_p}} (x_i) = n,  \mbox{ and } \ord_{ e_{D_p}} (y_i)  = m. 
\]
By definition the valuation $\nu_{D_p}$ is composed with the monomial 
valuation $\ord_{ e_{D_p}}$. This means that if $\psi_i: S_i \to S$ is the composition of 
modifications factoring $\pi$ and appearing in the Algorithm \ref{alg:tores}   until 
the point  $o_i $ appears in  $S_i$, then 
\begin{equation} \label{eq: sf3}
\nu_{D_p} (h ) = \ord_{ e_{D_p}} (  h  \circ \psi_i   ), 
\end{equation}
for any $h \in \mathcal{O}$. If $A$ is the plane curve defined by $h=0$ on $S$ 
we get that formula \eqref{eq: sf2} follows from \eqref{eq: sf} and \eqref{eq: sf3}. 
Formula \eqref{eq: sf4} follows  from Remark \ref{van-val} by the same argument.
\medskip

Recall that 
$ \underline{{\l}}_i  \in M_{R_i, L_i}$ 
and 
$ e_{D_p} = n e_{R_i} + m e_{L_i} $ belongs to the lattice $N_{R_i, L_i}$ which is dual to  $M_{R_i, L_i}$.
We get that
\begin{equation} \label{for:prod}
\langle  e_{D_p} ,  \underline{{\l}}_i \rangle = \langle  n e_{R_i} + m e_{L_i} ,  \l_{R_i} R_i  +   L_i   
\rangle = n \lambda_{R_i} + m. 
\end{equation}

Let $\F$ be a regular fan  of the latice $N_{R_i, L_i}$ which subdivides 
the cone  $\sigma_0^{R_i, L_i}$ and contains the ray $\R_{\geq 0}  \,  e_{D_p}$. 
Consider the modification $\phi:= \psi_{R_i, L_i}^\F : S^\F_i \to S_i$ (see \eqref{toroidal-map}). 
There is a chart of $\phi$ of the form 
\[
x_i = z^n t^a, \quad y_i = z^m t^b, \mbox{ with } n b- a m = \pm 1,
\]
where 
the representing divisor $D_p$ is defined by $z=0$. 
We deduce from this that  the order of vanishing of the Jacobian of $\phi$ along $D_p$ is equal to $n + m -1$. 

By definition the order of vanishing of the Jacobian of $\psi_i$ along $R_i$ is equal to $\l_{R_i} -1$.
By the chain rule the order of vanishing of the Jacobian of 
$\psi_i \circ \phi$ along $D_p$ is equal to 
\[
(\l_{R_i} -1) \cdot n + n + m -1 = n \l_{R_i}  + m -1. 
\]
The log-discrepancy of $D_p$ is equal to  the order of vanishing of the Jacobian of $\psi_i \circ \phi$ along $D_p$ plus one, that is,  it is equal to $n \l_{R_i}  + m$. 
Then, formula \eqref{eq: sf5} follows  by \eqref{for:prod}.
\end{proof}

\subsection{Main results}
The following theorem is a  generalization of 
Howald's description of  multiplier ideals of functions which are nondegenerate with respect their  Newton 
polyhedra (see \cite{Howald:Non-degenerate}) or  his 
description of the multiplier ideals of monomial ideals
(see \cite{MMI}). 

\begin{theorem}\label{Thm:Howalditerated} Let $C = \sum_{i=1}^r a_i C_i$, $a_i \in \N^*$ be a plane curve singularity at a point $o$ of  a smooth surface $S$. 
Let $R$ be a smooth branch which is not a component of $C$ and $\pi: \Sigma \to S $
a  toroidal embedded 
resolution of $C$  given by Proposition 
\ref{Prop: toroidal res}.  
Then, for any rational number $\xi >0$, we have  
\begin{equation} \label{eq: Ninc}
\mathcal{J}(\xi C)_o  = 
	\{ h \in \mathcal{O} \mid 
	\Newton_{R_i, L_i} (C_h) +  \underline{\l}_i \subset \mathrm{Int} (\xi \Newton_{R_i, L_i} (C)), \mbox{ for } i \in I_\pi
	 \} .
\end{equation}
\end{theorem}
\begin{proof}
We use 
\eqref{eq:form MI} to obtain that:
\begin{align} \label{form-4}
\mathcal{J}(\xi C)_o = 
	\{ h \in \mathcal{O} \mid  \nu_{D} (C_h)  + \l_{D} > \xi \nu_{D}(C)
	\mbox{ for } D  \in {E} (\pi)  \cup \{ C_1, \dots, C_r \}  \, \}.
\end{align}
We translate condition \ref{form-4} in terms of the fan tree (see Remark \ref{rem:tree}).
Let us take a segment
$[e_{R_i}, e_{L_i}]$ for some $i \in I_\pi$  of the decomposition 
$\mathcal{D}_\pi (C)$ of the fan tree $\theta_\pi (C)$  (see Definition \ref{def:fantreetr}).
By \eqref{eq: sf2}, for any marked point $p$ of $\theta_\pi (C)$ 
lying on the segment $[e_{R_i}, e_{L_i}]$ we have 
\begin{equation}\label{eq:mpbis}
\nu_{{D_p}} (C_h)  = \Phi_{\Newton_{R_i, L_i} (C_h)} ( e_{D_p}).
\end{equation}
In particular, if $C = C_h$ we  get that 
\begin{equation}
\xi \nu_{{D_p}}(C) = 
\xi \Phi_{\Newton_{R_i, L_i} (C)} ( e_{D_p}) \stackrel{\eqref{homothetic}}{=}  \Phi_{\xi \Newton_{R_i, L_i} (C)} ( e_{D_p}).
\end{equation}
By \eqref{eq: sf5}  and \eqref{eq:mpbis} we get
\begin{equation}
\nu_{{D_p}} (C_h)  + \lambda_{{D_p}}  = 
\Phi_{\Newton_{R_i, L_i} (C_h)} ( e_{D_p}) 
+ \langle  e_{D_p}, \underline{\l}_i \rangle = 
 \Phi_{\Newton_{R_i, L_i} (C_h) +  \underline{\l}_i } ( e_{D_p}).
\end{equation}
It follows that the condition:
\[
\nu_{{D_p}} (C_h)  + \lambda_{{D_p}}  > \xi \nu_{{D_p}}(C), 
\]
for any marked point $p$ in the segment $[e_{R_i}, e_{L_i}]$,
 is equivalent to: 
\begin{equation} \label{eq: cond2}
 \Phi_{\Newton_{R_i, L_i} (C_h) +  \underline{\l}_i } ( e_{D_p}) > \Phi_{\xi \Newton_{R_i, L_i} (C)} ( e_{D_p}).
\end{equation}
By Lemma \ref{lem :inclusion},   the expression \eqref{eq: cond2}  for $p \in [e_{R_i}, e_{L_i}]$,
is equivalent to the inclusion: 
\begin{equation} \label{eq: cond3}
\Newton_{R_i, L_i} (C_h) +  \underline{\l}_i \subset \mathrm{Int} (\xi \Newton_{R_i, L_i} (C)).
\end{equation}
Taking this into account  for every $i \in I_\pi$ 
ends the proof of \eqref{eq: Ninc}.
\end{proof}

\begin{cor} \label{prop-form}
 With the hypothesis and notation of Theorem \ref{Thm:Howalditerated} we have that:
\begin{align} \label{form-2}
\mathcal{J}(\xi C)_o = 
	\{ h \in \mathcal{O} \mid  \nu_{D} (C_h)  + \l_{D} > \xi \nu_{D}(C)
		\mbox{ for } D  \in \mathcal{R}_\pi (C) \cup \{ C_1, \dots, C_r \}  \, \}.	
\end{align}
In addition, if $C$ is reduced and $0 < \xi < 1$ we have
\begin{align} \label{form-3}
\mathcal{J}(\xi C)_o = 
	\{ h \in \mathcal{O} \mid  
	\nu_{D} (C_h)  + \lambda_{D}  > \xi \nu_{D}(C) \mbox{ ,  for } D \in \mathcal{R}_\pi (C)  \}.
	\end{align}
\end{cor}
\begin{proof} 
We keep the notation of the proof of Theorem \ref{Thm:Howalditerated} assuming that 
 $\pi: \Sigma \rightarrow S$ is  a maximal contact toroidal embedded 
resolution of $C$.

Notice that if $p =  e_{L_i}$ and $L_i$  is not a branch of $C$
then condition \eqref{eq: cond2} is always satisfied since $\l_{L_i} = 1$ and $\nu_{L_i} (C) =0$.
The same happens if $p = e_{R_1}$ since $R$ is not a component of $C$ by hypothesis. 
By Lemma \ref{lem :inclusion}, the inclusion \eqref{eq: cond3} is equivalent the inequality 
\eqref{eq: cond2} for  $p$ running through the ramification points 
of the tree $\theta_\pi (C)$ which belong to the segment 
$[e_{R_i}, e_{L_i}]$.
This proves formula \eqref{form-2}.
\medskip 

Assume now that $C$ is reduced, that is, $a_j =1$ for $j=1, \dots, r$. 
If $0 < \xi < 1$ and if $p = e_{L_i}$ is a branch of $C$ 
then condition \eqref{eq: cond2} always holds since $\l_{L_i} = 1$ and $\nu_{L_i}(C)  =1$. 
This implies \eqref{form-3}.
\end{proof}

We describe now the jumping numbers and the generators of the multiplier ideals: 
\begin{theorem}  \label{Thm:Monomiality} 
Let $C = \sum_{i=1}^r a_i C_i$, $a_i \in \N^*$ be a plane curve singularity at a point $o$ of  a smooth surface $S$. 
Let $R$ be a smooth branch which is not a component of $C$ and $\pi: \Sigma \to S $
a  toroidal embedded 
resolution of $C$  given by Proposition 
\ref{Prop: toroidal res}.  
 Denote by $\{ x_0, \dots, x_s \}$ the set obtained by taking a defining function for 
 every irreducible component of the completion $\hat{C}_\pi$.
We associate with a monomial $\mathcal{M}$ in $x_0, \dots, x_s$
the number
\begin{equation} \label{eq:mon}
\xi_{\mathcal{M}}:= 
		{\min} 
		\bigg\{ \dfrac{\nu_{D}(\mathcal{M})+\lambda_{D}}{\nu_{D}(C)}  \mid D \in \mathcal{R}_\pi (C) \cup \{ C_1, \dots, C_r \}
		\bigg\}.
\end{equation}
Then: 
\begin{enumerate}
\item
For any rational number $\xi >0$,  the multiplier ideal  $\mathcal{J}(\xi C)_o$  is generated by the finite set of monomials $\mathcal{M}$
in $x_0, \dots, x_s$ such that 
$\xi < \xi_\M \leq \xi+ 1$.

\item
The jumping numbers of 
the multiplier ideals of $C$ is the set of rational numbers 
$\xi_{\mathcal{M}}$ for $\mathcal{M}$ running through the monomials in $x_0, \dots, x_s$.
\end{enumerate}
\end{theorem}
\begin{proof}  
By  Corollary \ref{prop-form}
we have that the multiplier ideal  $\mathcal{J} ( \xi C)_o$
is a valuation ideal with respect to
the divisorial valuations $
\nu_{D}$, for $ D \in \mathcal{R}_\pi (C)$ 
and the vanishing order valuations
$\nu_{C_j}$, for $j = 1, \dots, r$.
By Remark \ref{rem:tree}, the set of components of 
   $\hat{C}_\pi$ contains a set of maximal contact curves of the pair $(\pi, C)$.
Corollary \ref{Cor:DGN} implies that 
that $\mathcal{J} ( \xi C)_o$ is generated by 
monomials in $x_0, \dots, x_s$.
By \eqref{form-4}
a monomial $\mathcal{M}$ in $x_0, \dots, x_s$
belongs to the multiplier ideal $\mathcal{J} ( \xi C)_o$ 
if and only if $\xi < \xi_{\mathcal{M}}$. 
\medskip 

Assume that  $\xi +1 < \xi_{\mathcal{M}} $ and let us prove that $\mathcal{M}$ is not a generator of the ideal $\mathcal{J} ( \xi C)_o$.  Our assumption implies that $1 < \xi_{\mathcal{M}}$ 
therefore 
$\mathcal{M}$ belongs to $\mathcal{J} (1\cdot C) = (f_C)$,
where  $f_C$ is a defining function of $C$ (see Lemma \ref{Lem:periodicidad}).
By Definition \ref{mcc}, we can assume that $f_C$
is a  monomial in 
$\{ x_0, \dots, x_s \}$.
It follows that 
$\mathcal{M} = \mathcal{M'}  \cdot f_C$ where $\mathcal{M'}$ is a monomial in $x_0, \dots, x_s$. 
Then, we check that $\xi_{\mathcal M} = 1 + \xi_{\mathcal{M}'}$. 
We iterate this argument 
to obtain
a monomial $\M''$ in 
$\{ x_0, \dots, x_s \}$
dividing  $\M$
and such that 
$\M''
\in \mathcal{J} (\xi C)_o$,   and
$\xi < \xi_{\M''} \leq \xi +1 $. This ends the proof of the first assertion.

\medskip 

Let us prove the second assertion. If $\xi >0$ is a jumping number of the multiplier ideals of $C$ there exists   $0 < \epsilon$ small enough,  
such that 
we have 
\[
\mathcal{J} (\xi C )_{o} \subsetneq \mathcal{J} (\xi' C )_{o}  \mbox{ and }  \mathcal{J} (\xi' C )_{o}  = \mathcal{J} (\xi'' C )_{o}, 
\]
for any $\xi', \xi'' \in (\xi - \epsilon, \xi)$.
There exists a monomial $\mathcal{M}$ in  $x_0, \dots, x_s$
such that $\mathcal{M} \in \mathcal{J} (\xi' C )_{o}$ for $\xi'  \in (\xi - \epsilon, \xi)$
and $\mathcal{M} \notin \mathcal{J} (\xi C )_{o}$. 
By   \eqref{form-2} 
for any 
 $D \in \mathcal{R}_\pi (C) \cup \{ C_1, \dots, C_r \} $
we have 
\[
\nu_{D} (\mathcal{M})+\lambda_{D} > \xi' \nu_{D} (C), \mbox{ for } \xi'  \in (\xi - \epsilon, \xi),
\]
while there exists  $D_0 \in \mathcal{R}_\pi (C) \cup \{ C_1, \dots, C_r \} $ 
such that the condition 
$
\nu_{D_{0}} (\mathcal{M})+\lambda_{{D_{0}}} > \xi \nu_{D_{0}}(C)
$
is not satisfied. 
This shows that $\xi = \xi_{\mathcal M}$ by continuity. 
Conversely, 
for every monomial $\M$ the number $\xi_{\M}$ is 
a jumping number since $\mathcal{M} \in \mathcal{J} (\xi' C )_{o}$ for $\xi'  \in (\xi_\M - \epsilon, \xi_\M)$
for $0 < \epsilon$ small enough, while $\mathcal{M} \notin \mathcal{J} (\xi_{\M} C )_{o}$. 
\end{proof}

\begin{remark} \label{rem:red}
The toroidal embedded resolution $\pi$ is the composition of $\ell_\pi$ 
regularized Newton modifications. 
Assume that $\pi$ is such that the integer $\ell_\pi$ is 
the smallest one. 
Then, by  \cite[Th. 3.12]{LO} one has that 
$\pi$
 is the minimal 
embedded resolution of $C$ and
the set 
$\{ x_0, \dots, x_s \}$ in Theorem 
\ref{Thm:Monomiality} 
is obtained by taking a defining function for every element in a sequence of maximal contact curves of the pair $(\pi, C)$. 
Up to relabelling, let us denote by 
$\{ x_0, \dots, x_{s'} \}$ 
a subset $\{ x_0, \dots, x_s \}$ consisting of the defining functions of a set of maximal contact curves of $\pi$.
If  $C$ is reduced 
and $0 < \xi < 1$ then the multiplier ideal $\mathcal{J}(\xi C)_o$  is generated by monomials in the sequence 
$\{ x_0, \dots, x_{s'} \}$. 
This is consequence of \eqref{form-3} and Corollary  \ref{Cor:DGN}.
Then, reasoning as in the proof of Theorem \ref{Thm:Monomiality} we obtain 
that
if $0 < \xi_0 < 1$ 
is a jumping number there 
exists a monomial $\M$ in 
$x_0, \dots, x_{s'}$
such that $\xi_0$ equals 
\begin{equation} \label{xiMred}
{ \xi_{\mathcal{M}}^{\mathrm{red}} }:= 
	{\mathrm{min}} 
	\{ ({\nu_{D}(\mathcal{M})+\lambda_{D}}) ({\nu_{D}(C)})^{-1}  \mid D \in \mathcal{R}_\pi (C) 
	\}.
\end{equation}
Conversely, if $\M$ is a monomial in 
$x_0, \dots, x_{s'}$ such that  $0 < \xi_{\M}^{\mathrm{red}}< 1$ then $\xi_{\M}^{\mathrm{red}}$ is a jumping number of $C$. 
Notice that if 
$\xi'> 1$ is not an integer then $0 < \xi:= \xi' - \lfloor \xi' \rfloor <1$ and by Lemma 
\ref{Lem:periodicidad}, we have that 
$\mathcal{J}(\xi' C)_o = f^{\lfloor \xi' \rfloor } \cdot  \mathcal{J}(\xi  C)_o$.
\end{remark}

\begin{remark} 
A different proof of the monomiality in Theorem \ref{Thm:Monomiality} 
was given in \cite{Blanco}. Their proof holds more generally 
for integrally closed ideals of $\mathcal{O}$ and  uses  
the correspondence between antinef  divisors and complete ideals. The approach of  this paper may be also generalized for the study of the multiplier ideals of ideals of $\mathcal{O}$ (see  \cite[Sect. 5.3]{PhD-Robredo}).

\end{remark}

\begin{remark}
\label{rem:char}
The results of this paper  also hold for an algebroid curve $(C,o)$ on a smooth surface $S$ over algebraically closed field of 
arbitrary characteristic. Notice that in this case we have a unique minimal embedded resolution of a plane curve $C$ 
and any other embedded resolution factors through it. The argument given in \cite[Th. 9.2.18]{PAG} implies 
that the multiplier ideals $\mathcal{J} (\xi C)_o$  are independent of the choice of embedded resolution. 
Then, we use that the toroidal embedded resolutions 
can be built in this setting independently of the characteristic of the base field. 
In particular, there is no need to pass through Newton-Puiseux series 
in order to connect the fan tree $\theta_\pi (C)$ 
of a toroidal embedded resolution with the functions in the valuative tree 
(see \cite[Sec. 1.6.6]{GBGPPP20} and \cite{GGPValTree} for details).
\end{remark}

\begin{remark}
The slope funtion of $\slp_\pi$ of the fan tree $\theta_\pi (C)$ was introduced in  \cite{GBGPPP20}, see Definition \ref{def:fantreetr}. This function 
determines explicitly the values of the log-discrepancies of 
a rupture component $R$ of $\pi^*(C)$ (see \cite[Prop. 8.16 (2)]{GGPValTree}) and the values of 
the divisorial valuation $\ord_{R}$ at the branches of the completion $\hat{C}_\pi$ (see 
\cite[Cor. 3.26, Prop. 7.18 and Prop. 8.16 (1)]{GGPValTree}).
That is all the data required to compute the multiplier ideals of $C$,
as in Example \ref{Ex:2Branches} below. 
In \cite[Sec. 6.5]{GBGPPP20} 
it is shown 
how to identify 
the fan tree $\theta_{\pi} (C)$ with
the Eggers-Wall tree 
of $\hat{C}_\pi$.
This identification is related to the embedding of the fan tree 
$\theta_{\pi} (C)$ in the \textit{valuative tree} 
$\mathcal{V}_{R}$ 
of normalized semivaluations 
with respect to  the smooth branch $R$ (see \cite{GGPValTree}), which 
sends 
a rational point $p \in \theta_\pi (C)$ to the normalized valuation $\frac{\nu_{D_p}}{\nu_{D_p} (R)}$. 
The valuative tree $\mathcal{V}_{R}$ has been intensively studied by Favre and Jonsson in \cite{TVT}, see also \cite[Section 7]{JonssonDynamicsBerkovich}.
\end{remark}

\begin{example} \label{Ex:2Branches}
Let us consider the toroidal embedded resolution $\pi$ of the plane curve  $C$ of Example \ref{ex:2branches}. 
The set of rupture components of $\pi^*(C)$ is $\mathcal{R}_\pi (C) = \{R_1, R_2, R_3 \}$.
The log-discrepancies of the exceptional divisors $R_1$, $R_2$ and $R_3$ are 
 \begin{equation} \label{eq:lambda}
\lambda_{R_2} = 5, 
 \lambda_{R_3} = 8,  \mbox{ and }\lambda_{R_4} = 13. 
 \end{equation}
Table \ref{Table:ValuesEx:2Branches} provides the
 the values of the divisorial valuations $\nu_{R_i}$  of the rupture components $R_i$, for $i=2,3,4$,  
at  the branches of $\hat{C}_\pi$.

\begin{table}[h!]
 \begin{center}
\begin{tabular}{c|c|c|c|c|c|c}
                              &   $R_1$ & $L_1$ & $L_2$  & $C_{1}$ & $C_{2}$ & $C$  \\  \hline
{$\nu_{R_{2}}$} &   $2$ & $3$ & $6$  & $12$    & $9$     & $21$               \\  
\hline
{$\nu_{R_{3}}$} &   $3$ & $5$ & $9$  & $18$    & $15$    & $33$                \\ 
\hline
{$\nu_{R_{4}}$} &   $4$ & $6$ & $15$ & $30$    & $18$    & $48$               \\ 
\end{tabular}
\caption{List of values for the divisorial valuations in Example \ref{Ex:2Branches}.}
\label{Table:ValuesEx:2Branches}
 \end{center}
\end{table}
 We illustrate in this example how the computation of the set of jumping numbers smaller than one 
and of a system of generators of the corresponding 
multiplier ideals reduces to  an optimization problem 
in terms of the data of values of Table \ref{Table:ValuesEx:2Branches} and the log-discrepancies \eqref{eq:lambda}. 

\medskip

Recall from Example \ref{Ex:2BranchesFanTree} that
$R_1 = C_x$, $L_1 = C_y$ and $L_2 = C_{z}$, where $z = y^2 + x^3$. 
By Theorem \ref{Thm:DGN}, the functions $x$, $y$ and $z$ define a minimal generating sequence of 
the divisorial valuations $\nu_{R_2}$, $\nu_{R_3}$ and $\nu_{R_4}$. 
By Corollary \ref{Thm:Monomiality}  and Remark \ref{rem:red}, 
if  $0<\xi<1$, if $\xi$ is a jumping number then
there exists a monomial $\M = x^{a}y^{b}z^{c}$
such that $\xi = \xi^{red}_\M$ where 
\[ 
\xi^{red}_{\M} := \min \{ {(2a+3b+6c+5)}/{21}, {(3a+5b+9c+8)}/{33},{(4a+6b+15c+13)}/ {48} \}.
\]
The multiplier ideal
$\mathcal{J}(\xi C)_o$, 
is generated by monomials in $x$, $y$ and $z$
such that $\xi < \xi^{red}_\M \leq \xi +1$
We indicate below the multiplier ideals ideals $\mathcal{J}(\xi C)_o$ for $\xi$ running through the 30
jumping numbers in the interval $(0,1)$, where we have underlined a monomial $\M= x^{a}y^{b}z^{c} \in \mathcal{J}(\xi C)_o$
 if  $\xi^{red}_\M > \xi$ is the next jumping number after $\xi$.
The output is not a minimal set of generators of the ideals, for instance, we get
$\mathcal{J}\left( \frac{5}{21} C\right)_o = (x, y, z)$, which is equal $(x, y)$.
\end{example}

\[
\begin{array} {lcl}
\mathcal{J} \left(\frac{5}{21} \cdot C\right)_o & = &  (\underline{{x}}, y, z)
\\
\mathcal{J} \left( \frac{1}{3}  \cdot C\right)_o & = & (x^2,\underline{{y}}, z)
\\
\mathcal{J} \left( \frac{8}{21}  \cdot C\right)_o & = & (\underline{{x^2}},xy,y^2, z)
\\
\mathcal{J} \left(  \frac{14}{33}  \cdot C\right)_o & = & (x^3,\underline{{xy}},y^2, z)
\\
\mathcal{J} \left( \frac{10}{21}  \cdot C\right)_o & = & (\underline{{x^3}},x^2y,y^2, \underline{{z}})
\\
\mathcal{J}\left( \frac{17}{33} \cdot C \right)_o & = & (x^4, x^2y, \underline{{y^2}}, xz, yz, z^2)\\
\mathcal{J}\left( \frac{25}{48} \cdot C \right)_o & = & (x^4, \underline{{x^2y}}, xy^2, y^3, xz, yz, z^2)\\
\mathcal{J}\left( \frac{27}{48} \cdot C \right)_o & = & (\underline{{x^4}}, x^3y, \underline{{xy^2}}, y^3, xz, yz, z^2)\\
\mathcal{J}\left( \frac{29}{48} \cdot C \right)_o & = & (x^5, x^3y, x^2y^2, y^3, \underline{{xz}}, yz, z^2)\\
\mathcal{J}\left( \frac{20}{33} \cdot C \right)_o & = & (x^5, \underline{{x^3y}}, x^2y^2, \underline{{y^3}}, x^2z, yz, z^2)\\
\mathcal{J}\left( \frac{31}{48} \cdot C \right)_o & = & (x^5, x^4y, x^2y^2, xy^3, y^4, x^2z, \underline{{yz}},  z^2)\\
\mathcal{J}\left( \frac{2}{3} \cdot C \right)_o & = & (\underline{{x^5}}, x^4y, \underline{{x^2y^2}}, xy^3, y^4, x^2z, xyz, y^2z, z^2)\\
\mathcal{J}\left( \frac{33}{48} \cdot C \right)_o & = & (x^6, x^4y, x^3y^2, xy^3, y^4, \underline{{x^2z}}, xyz, y^2z, z^2)\\
\mathcal{J}\left( \frac{23}{33 } \cdot C \right)_o & = & (x^6, \underline{{x^4y}}, x^3y^2, \underline{{xy^3}}, y^4, x^3z, xyz, y^2z,  z^2)\\
\mathcal{J}\left( \frac{35}{48} \cdot C \right)_o & = & (x^6, x^5y, x^3y^2, x^2y^3, y^4, x^3z, \underline{{xyz}}, y^2z, z^2)\\
\mathcal{J}\left( \frac{25}{33} \cdot C \right)_o & = & (\underline{{x^6}}, x^5y, \underline{{x^3y^2}}, x^2y^3, \underline{{y^4}}, x^3z, x^2yz, y^2z, z^2)\\
\mathcal{J}\left( \frac{37}{48} \cdot C \right)_o & = & (x^7, x^5y, x^4y^2, x^2y^3, xy^4,y^5, \underline{{x^3z}}, x^2yz, y^2z, \underline{{z^2}})\\
\mathcal{J}\left( \frac{26}{33} \cdot C \right)_o & = & (x^7, x^5y, x^4y^2, x^2y^3, xy^4,y^5, x^4z, x^2yz, \underline{{y^2z}}, xz^2, yz^2)\\
\mathcal{J}\left( \frac{17}{21} \cdot C \right)_o & = & (x^7, \underline{{x^5y}}, x^4y^2, \underline{{x^2y^3}}, xy^4, y^5, x^4z, x^2yz, xy^2z,y^3z, xz^2, yz^2)\\
\mathcal{J}\left( \frac{39}{48 } \cdot C \right)_o & = & (x^7, x^6y, x^4y^2, x^3y^3, xy^4, y^5, x^4z, \underline{{x^2yz}}, xy^2z, y^3z, xz^2, yz^2)\\
\mathcal{J}\left( \frac{28}{33 } \cdot C \right)_o & = & (\underline{{x^7}}, x^6y, \underline{{x^4y^2}}, x^3y^3, \underline{{xy^4}}, y^5, x^4z, x^3yz, xy^2z,y^3z, xz^2, yz^2)\\
\mathcal{J}\left( \frac{41}{48} \cdot C \right)_o & = & (x^8, x^6y, x^5y^2, x^3y^3, x^2y^4, y^5, \underline{{x^4z}}, x^3yz, xy^2z,y^3z, \underline{{xz^2}}, yz^2)\\
\mathcal{J}\left( \frac{29}{33 } \cdot C \right)_o & = & (x^8, \underline{{x^6y}}, x^5y^2, \underline{{x^3y^3}}, x^2y^4, \underline{{y^5}}, x^5z, x^3yz, xy^2z,y^3z, x^2z^2, yz^2)\\
\mathcal{J}\left( \frac{43}{48 } \cdot C \right)_o & = & (x^8, x^7y, x^5y^2, x^4y^3, x^2y^4, xy^5, y^6, x^5z, x^3yz, \underline{{xy^2z}}, y^3z,  x^2z^2, yz^2)\\
\mathcal{J}\left( \frac{19}{21 } \cdot C \right)_o & = & (\underline{{x^8}}, x^7y, \underline{{x^5y^2}}, x^4y^3, \underline{{x^2y^4}}, xy^5, y^6, x^5z, x^3yz, x^2y^2z,y^3z, x^2z^2, yz^2)\\
\mathcal{J}\left( \frac{45}{48 } \cdot C \right)_o & = & (x^9, x^7y, x^6y^2, x^4y^3, x^3y^4, xy^5, y^6, x^5z, \underline{{x^3yz}}, x^2y^2z,y^3z, x^2z^2, \underline{{yz^2}})\\
\mathcal{J}\left( \frac{31}{33 } \cdot C \right)_o & = & (x^9, x^7y, x^6y^2, x^4y^3, x^3y^4, xy^5, y^6, x^5z, x^4yz, x^2y^2z, \underline{{y^3z}}, x^2z^2, xyz^2, y^2z^2)\\
\mathcal{J}\left( \frac{20}{21} \cdot C \right)_o & = & (x^9, x^7y, x^6y^2, x^4y^3, x^3y^4, xy^5, y^6, \underline{{x^5z}}, x^4yz, x^2y^2z, xy^3z, y^4z, \underline{{x^2z^2}}, xyz^2, y^2z^2)\\
\mathcal{J}\left(  \frac{32}{33} \cdot C\right)_o& = & (x^9, \underline{{x^7y}}, x^6y^2, \underline{{x^4y^3}}, x^3y^4, \underline{{xy^5}}, y^6, x^6z, x^4yz, x^2y^2z, xy^3z, y^4z, x^3z^2, xyz^2, y^2z^2)\\
\mathcal{J}\left(  \frac{47}{48} \cdot C \right)_o & = & (x^9, x^8y, x^6y^2, x^5y^3, x^3y^4, x^2y^5, y^6, x^6z, x^4yz, x^2y^2z, xy^3z, y^4z, x^3z^2, xyz^2, y^2z^2)\\
\end{array}
\]

\medskip

\begin{remark} \label{Python}
The fourth author wrote a program in Python to compute 
the jumping numbers and finite presentations for the multiplier ideals 
of a plane curve singularity 
in terms the log-discrepancies of the exceptional divisors $D \in \mathcal{R}_\pi (C)$, and 
the values of $\nu_{D}$ on the components of $\hat{C}_\pi$. 
The code is available at 
\href{https://github.com/tdimitch/jumping-numbers}{https://github.com/tdimitch/jumping-numbers}.
\end{remark}

\begin{ack}
We are grateful to Patrick Popescu-Pampu for his comments on  a preliminary version of this paper. 
We would like to thank the referee for the suggestions which have helped us to improve the presentation of the paper.
\end{ack}

\providecommand{\href}[2]{#2}


\begin{thebibliography}{99}

\bibitem{AAD}
{Alberich-Carrami\~{n}ana}, M.,  { \'{A}lvarez Montaner}, J.,  {Dachs-Cadefau}, F. 
(2016). 
{Multiplier ideals in two-dimensional local rings with rational singularities}. 
  \textit{Michigan Math. J.} 
{65}(2):287--320. 
  DOI: 10.1307/mmj/1465329014

\bibitem{Blanco}
{Alberich-Carrami\~{n}ana}, M.,  { \'{A}lvarez Montaner}, J., {Blanco}, G. 
(2021).
 {Monomial generators of complete planar ideals}.
 \textit{J. Algebra Appl.}
  20(3):Article 2150032.   DOI: 10.1142/S0219498821500328
 
\bibitem{ACO96}
{A'Campo}, N.,  
{Oka}, M.  
(1996).
{Geometry of plane curves via {T}schirnhausen resolution tower}.
\textit{Osaka J. Math.}
 {33} (4): 1003--1033. 
 
\bibitem{Berkesch10}
{Berkesch}, C., 
{Leykin}, A.  
(2010).
{Algorithms for Bernstein--Sato polynomials and multiplier ideals}. 
\textit{Proceedings of the 2010 International
  Symposium on Symbolic and Algebraic Computation}. 
  New York, NY:
   Association for Computing Machinery, 
  pp. 99--106.

\bibitem{BudurSingulInvMilnorFiber}
{Budur},   N. 
(2012).
{Singularity invariants related to {M}ilnor fibers: survey},
In: 
Campillo, A.,
Cardona, G.,
Melle-Hern\'andez, A.,
Veys, W.,
 Zúñiga-Galindo, A.
 ed(s).
  \textit{Zeta functions in algebra and geometry}. 
  Contemp. Math. 566.  Providence, RI: American Mathematical Society, pp.161--187.

\bibitem{CG03}
{Campillo}. A., 
{Galindo}, C. 
(2003).
{The {P}oincar\'{e} series associated with finitely many monomial valuations}.
  \textit{Math. Proc. Cambridge Philos. Soc.}
  {134}(3):433--443. 
  DOI: 0.1017/S030500410200645X.

\bibitem{Cox}
Cox, D.~A., 
Little, J.~B.,
Schenck,  H.~K. 
(2011). 
 {Toric varieties}, 
  Graduate Studies in Mathematics, vol. 124.
 Providence, RI: American Mathematical Society.

\bibitem{CNL14}
{Cassou-Nogu{\`e}s}, Pi.,  {Libgober}, A. 
(2014).
{Multivariable {H}odge theoretical invariants of germs of plane curves. {II}}, 
  In: 
Campillo, A., 
Kuhlmann, F.-V., 
Teissier, B. 
ed(s).
\textit{Valuation Theory in Interaction}.
 EMS Ser. Congr. Rep. 
 Z\"{u}rich: Eur. Math. Soc.,
  pp.~82--135. 

\bibitem{DelGalNuGenSeq}
Delgado, F., Galindo, C., {N\'{u}\~{n}ez}, A. (2008)
{Generating sequences and {P}oincar\'{e} series for a finite set of plane divisorial valuations}. 
\textit{Adv. Math.}
219(5):1632--1655. DOI: 10.1016/j.aim.2008.06.017.



\bibitem{ResolComplex} {D\~ung Tr\'ang}. L., Oka, M.  (1995).
{On resolution complexity of plane curves}, 
\textit{Kodai Math. J.} 18 (1):1--36. 

\bibitem{ELSVJumpingCoefficients}
Ein, L., 
Lazarsfeld, R.,  
Smith, K.~E.  
Varolin, D.
(2004).
 {Jumping coefficients of multiplier ideals}. 
 \textit{Duke Math. J.}
  123(3):469--506.


\bibitem{EwCombConv}
Ewald, G.
(1996).
\textit{Combinatorial convexity and algebraic geometry}, 
New York, NY: Springer-Verlag.

\bibitem{TVT}
Favre, C., 
Jonsson,  M. 
(2004). 
\textit{The valuative tree}, Lecture Notes in Mathematics, vol. 1853.
Berlin: Springer-Verlag.

\bibitem{FJ05}
Favre, C., 
Jonsson,  M. 
(2005). 
 {Valuations and multiplier ideals}, 
 \textit{J. Amer. Math. Soc.} 
 18 (3):655--684. 
 
\bibitem{FTV}
Fulton, W.  (2016). 
\textit{Introduction to toric varieties}
Annals of Mathematics Studies, vol. 131.
Princeton: Princeton University Press.



\bibitem{GuzmanPhD}
Guzm{\'a}n~Dur{\'a}n, C.R. (2018).
 {Ideales multiplicadores de curvas planas irreducibles}.
  Ph.D. dissertation.
   Centro de investigaci{\'o}n en matem{\'a}ticas, Guanajuato, Mexico.
   Available at: \url{http://cimat.repositorioinstitucional.mx/jspui/handle/1008/725}

\bibitem{GGPValTree}
{Garc{\'\i}a Barroso}, E.R., 
{Gonz{\'a}lez Perez}, P.D.,
  {Popescu-Pampu}, P. 
  (2019).
   {The valuative tree is the projective limit of
  {E}ggers-{W}all trees}.
  \textit{Rev. R. Acad. Cienc. Exactas F\'{\i}s. Nat. Ser. A Mat. RACSAM}.
 113(4):4051--4105. 
 
\bibitem{GBGPPP20}
{Garc{\'\i}a Barroso}, E.R., 
{Gonz{\'a}lez Perez}, P.D.,
  {Popescu-Pampu}, P. 
(2020).
   {The combinatorics of plane curve singularities. How Newton polygons blossom into lotuses}.
  In: Cisneros-Molina, J.~L., Dung Tr\'{a}ng, L., Seade, J. ed(s).
  \textit{Handbook of Geometry and
Topology of Singularities I}.
Cham: Springer, pp.~1--150.

\bibitem{GPRQo}
{Gonz\'alez P\'erez},  P.D. 
  (2003).
{Toric embedded resolutions of quasi-ordinary hypersurface singularities}.
 \textit{Ann. Inst. Fourier (Grenoble)}.
{53}(6):1819--1881. 

\bibitem{HJHNORdTree}
Hyry, E.,
J\"{a}rvilehto, T.~  
 (2011).
{Jumping numbers and ordered tree structures on the dual graph}. 
\textit{Manuscripta Math.} 
136(3):411--437.

\bibitem{HJFormulaJN}
Hyry, E.,
J\"{a}rvilehto, T.~ 
(2018).
{A formula for jumping numbers in a two-dimensional regular local ring}.
\textit{J. Algebra}
516:437--470. 

\bibitem{MMI}
Howald,  J.~A. 
(2001).
{Multiplier ideals of monomial ideals}, 
\textit{Trans. Amer. Math. Soc.} 
353(7):2665--2671.


\bibitem{Howald:Non-degenerate}
Howald,  J.~A. 
(2003).
{Multiplier Ideals of Sufficiently General Polynomials}, 
 arXiv:math/0303203 [math.AG].
 
\bibitem{JarJNSCI}
J{\"a}rvilehto,  T.~
(2011).
{Jumping numbers of a simple complete ideal in a two-dimensional regular local ring}, 
\textit{Mem. Amer. Math. Soc.}
214 (1009), pp:78. 

\bibitem{JonssonDynamicsBerkovich}
Jonsson, M.~  
(2015).
{Dynamics on Berkovich spaces in low dimensions}.
In: Ducros, A.,  Favre, C.,  Nicaise, J. ed(s).
\textit{Berkovich spaces and applications}. 
Lecture Notes in Math., vol. 2119.
Cham: Springer, pp.~205--366.


\bibitem{PAG}
Lazarsfeld, R. 
(2004).
 \textit{Positivity in algebraic geometry. {II}}.
 Ergebnisse der
   Mathematik und ihrer Grenzgebiete. 3. Folge, vol.~49.
  Heidelberg: Springer Berlin.
  
\bibitem{JNN}
Naie, D.~  
(2009).
{Jumping numbers of a unibranch curve on a smooth surface}.
\textit{Manuscripta Math.} 
128(1):33--49

\bibitem{LO}
 L\^{e}, D.~T., Oka, M., (1995).
{On resolution complexity of plane curves}.
\textit{Kodai Math. J.} 18(1): 1--36





\bibitem{OdaCB}
Oda, T. 
(1988).
\textit{Convex bodies and algebraic geometry}.
Ergebnisse der Mathematik
 und ihrer Grenzgebiete (3),
 vol.~15. 
 Berlin: Springer-Verlag.


\bibitem{OkaToroidalRes}
Oka, M.~
(1996).
{Geometry of plane curves via toroidal resolution}, 
In:  Campillo L\'opez, A., Narv\'aez Macarro, L. ed(s).
\textit{Algebraic Geometry and Singularities}. 
 Basel: Birkh\"{a}user, pp.~95--121.
 

\bibitem{PP03}
{Popescu-Pampu}, P.~
(2003).
{Approximate roots}. 
In: Kuhlmann, F.-V., Kuhlmann, S., Marshall, M. ed(s).
\textit{Valuation theory and its  applications, {V}olume {II}}.
Fields Inst. Commun.,  vol.~33.
  Providence, RI: Amer. Math. Soc.,  pp.~285--321. 

\bibitem{PhD-Robredo}
Robredo~Buces, M.
(2019).
 {Invariants of singularities, generating sequences and toroidal structures}
  Ph.D. dissertation.
   Universidad Complutense de Madrid, Madrid, Spain.
   Available at:  \url{https://www.icmat.es/Thesis/2019/Tesis_Miguel_Robredo.pdf}

\bibitem{Shibuta11}
{Shibuta},T.
(2011). 
 {Algorithms for computing multiplier ideals}, 
  \textit{J. Pure  Appl. Algebra}
  215(12): 2829--2842.
   
\bibitem{VS}
Spivakovsky, M.~
 (1990).
{Valuations in function fields of surfaces}, 
\textit{Amer. J. Math.} 
112(1):107--156.



\bibitem{STIrrelExcDiv}
Smith, K.E.,  Thompson, H.M. 
(2007).
{Irrelevant exceptional divisors for curves on a smooth surface},
In: Corso, A., Migliore, J, Polini, C. ed(s).
 \textit{Algebra, geometry and their interactions}.
 Contemp. Math., vol. 448. 
 Providence, RI:  Amer. Math. Soc.,  pp.~245--254. 

\bibitem{SchTuck}
{Schwede}, K.,  
{Tucker}, K. 
(2012).
{A survey of test ideals}. 
In: Francisco, C., Klingler, L., Sather-Wagstaff, S., Vassilev, J.C. ed(s).
\textit{Progress in commutative algebra 2}. 
Berlin:  Walter de Gruyter, pp.~39--99.


\bibitem{TuckerJNAlgSurf}
Tucker, K. 
(2010). 
{Jumping numbers and multiplier ideals on algebraic surfaces},
Ph.D. dissertation. University of Michigan, Ann Arbor, MI: ProQuest LLC.

\bibitem{TuckerJNRatSing}
Tucker, K. 
(2010). 
{Jumping numbers on algebraic surfaces with rational singularities}, 
\textit{Trans. Amer. Math. Soc.}
{362}(6):3223--3241. 

\bibitem{Zhang}
{Zhang}, M. 
(2019).
 {Multiplier ideals of analytically irreducible plane curves},  
 arXiv:1907.06281v3 [math.AG].
 

\end{thebibliography}
\end{document}